\pdfoutput=1
\documentclass[preprint,12pt]{elsarticle}

\usepackage{graphicx}

\usepackage{bm}
\usepackage{amsmath}
\usepackage{amssymb}
\usepackage{amsthm}
\usepackage{color}
\usepackage{ulem}

\newcommand{\IGNORE}[1]{}
\newcommand{\ignore}[1]{}
\newcommand{\veps}{\varepsilon}
\newcommand{\opn}{\operatorname}

\newcommand{\jt}{\textstyle}

\newcommand{\der}[2]{\frac{\partial #1}{\partial #2}}

\newcommand{\erf}{\opn{erf}}
\renewcommand{\H}{\mathcal{H}}
\newcommand{\V}{\mathcal{V}}

\newcommand{\x}{\tilde{x}}

\newcommand{\rb}[2]{\raisebox{#1pt}{$#2$}}
\newcommand{\la}{\langle}
\newcommand{\ra}{\rangle}

\newcommand{\e}[1]{{(#1)}}

\journal{Journal of Computational Physics}

\begin{document}

\begin{frontmatter}

\title{Accurate Spectral Numerical Schemes for Kinetic Equations with Energy Diffusion}

\author[berk]{Jon Wilkening}
\author[courant]{Antoine Cerfon}
\author[maryland]{Matt Landreman}

\address[berk]{Department of Mathematics, University of California,
Berkeley}
\address[courant]{Courant Institute of Mathematical Sciences, New York University}
\address[maryland]{Institute for Research in Electronics and Applied Physics,
University of Maryland}

\begin{abstract}
  We examine the merits of using a family of polynomials that are
  orthogonal with respect to a non-classical weight function to
  discretize the speed variable in continuum kinetic calculations.  We
  consider a model one-dimensional partial differential equation
  describing energy diffusion in velocity space due to Fokker-Planck
  collisions. This relatively simple case allows us to compare the
  results of the projected dynamics with an expensive but highly
  accurate spectral transform approach. It also allows us to integrate
  in time exactly, and to focus entirely on the effectiveness of the
  discretization of the speed variable. We show that for a fixed
  number of modes or grid points, the non-classical polynomials can be
  many orders of magnitude more accurate than classical Hermite
  polynomials or finite-difference solvers for kinetic equations in
  plasma physics. We provide a detailed analysis of the difference in
  behavior and accuracy of the two families of polynomials.
  For the non-classical polynomials, if the
  initial condition is not smooth at the origin when interpreted as a
  three-dimensional radial function, the exact solution leaves the
  polynomial subspace for a time, but returns (up to roundoff
  accuracy) to the same point evolved to by the projected dynamics in
  that time.  By contrast, using classical polynomials, the exact
  solution differs significantly from the projected dynamics solution
  when it returns to the subspace. We also explore the connection
  between eigenfunctions of the projected evolution operator and
  (non-normalizable) eigenfunctions of the full evolution operator, as
  well as the effect of truncating the computational domain.
\end{abstract}


\end{frontmatter}

\section{Introduction}

First-principles based descriptions of transport processes in plasmas
require the solution of high-dimensional kinetic equations for the
phase-space distribution function
\cite{hazeltine1,hazeltine2}. 
Often, diffusion in velocity-space plays an important physical role,
as discussed in \cite{alex, abel}, and references therein.
For example, in kinetic turbulence, there is a cascade of energy in velocity-space
in addition to real space, so velocity diffusion cuts off the cascade
at small velocity scales. In this sense, velocity diffusion
plays a role similar to that of viscosity in conventional hydrodynamic fluid turbulence.
Just as viscosity is important in hydrodynamic turbulence no matter how large the Reynolds number,
velocity diffusion is important in kinetic turbulence no matter how small the collisionality.
Velocity diffusion is essential to dissipate injected energy and thereby
permit a statistically steady state.

Solving these kinetic equations numerically is
computationally intensive \cite{candy1,barnes1}, so an
important aspect of the theoretical effort is to find new optimized
discretization schemes. While high order accurate discretization
schemes for the spatial variables have been succesfully used for many
years, finding an ideal discretization method remains
particularly challenging for the discretization of velocity space in
situations involving Fokker-Planck collisions \cite{helander}. Since
the Fokker-Planck collision operator has terms involving first and
second order derivatives with respect to the velocity variables, the
discretization method must allow accurate differentiation. The scheme
must also allow accurate integration since physical quantities such as
the number density, the mean fluid velocity and the pressure depend on
velocity moments of the distribution function.

Recently, promising new approaches based on spectral and
pseudo-spectral representations have been investigated
\cite{bratanov,landreman}. It was shown in \cite{bratanov} that a
Hermite representation for the parallel velocity has advantages over
the more common finite difference schemes used in numerical
simulations. In \cite{landreman}, different representations for the
speed coordinate are explored. It is found that because the variable
has values in $[0,\infty)$ instead of the entire real axis, a
little-known family of polynomials (see \cite{shizgal,ball,ghiroldi} and references
therein) gives much better performance than finite difference schemes
and schemes based on classical orthogonal polynomials. High accuracy
is obtained on very coarse grids for both differentiation and
integration of Maxwellian-like functions, which are the functions of
interest in many applications of plasma physics \cite{landreman}.

The purpose of this paper is to explore the
suitability of the non-classical polynomials for
initial-value calculations of turbulent plasma transport in the presence of collisions
\cite{abel,barnes2}. 
To do so, we consider a model one-dimensional problem
describing energy diffusion due to Fokker-Planck collisions~\cite{vsck1}:
 
\begin{align}\label{eq:model}
  \der{U}{t} = \frac{1}{x^2}\der{}{x}\left[
    \Psi(x)x^2 e^{-x^2}\der{}{x}\Big(e^{x^2}U\Big)\right],
  \qquad (x>0,\;t>0),
\end{align}
where
\begin{align}\label{eq:psi:def}
  \Psi(x) = \frac{1}{2x^3}\left[\erf(x) - \frac{2}{\sqrt{\pi}}xe^{-x^2}\right], \qquad
  \erf(x) = \frac{2}{\sqrt{\pi}}\int_0^x e^{-s^2}\,ds.
\end{align}
Here, $\Psi$ differs from the usual Chandrasekhar function by an
additional factor of $1/x$.  Our choice to focus on this particular
model problem is motivated by the following characteristic
features. First, the right-hand side of \eqref{eq:model} corresponds
exactly to the speed variable part of the energy diffusion piece in
the linearized Landau-Fokker-Planck operator for same-species
collisions \cite{abel,barnes2}. Since many state-of-the-art plasma
turbulence codes (e.g.~\cite{candy1,barnes1}) use variants of the
linearized Landau-Fokker-Planck operator to describe collisions, the
results presented here are directly relevant to the computational
effort to simulate transport processes in plasmas.  Second, the
relative simplicity of (\ref{eq:model}) makes it possible to represent
the solution semi-analytically using a spectral transform method
\cite{vsck1}, which we can then use to study the properties and
accuracy of various discretization schemes. Third, \eqref{eq:model}
has several physically satisfying properties. Any well-behaved initial
distribution function relaxes to a Maxwellian distribution function
$U\propto e^{-x^2}$ as $t\rightarrow\infty$ (the ``H-theorem''). Also,
for all $t$, $\partial/\partial t\left(\int Ux^{2}dx\right)=0$,
i.e.~the number of particles is conserved.  The best time-dependent
numerical schemes are designed to satisfy these basic properties
exactly~\cite{barnes2}.

While we consider a single velocity dimension in \eqref{eq:model}, any
accurate numerical scheme for this equation is immediately applicable to simulations with more
velocity dimensions and more complete collision operators, such as the operator
in \cite{abel} or the linearized Fokker-Planck operator, for the following reasons.
For either of these collision operators in spherical velocity coordinates,
the only term involving any $\partial/\partial x$ derivatives of the
distribution function is the right-hand side of \eqref{eq:model}.
For example, in the pitch-angle diffusion term, the speed $x$ appears only as a parameter,
not as a derivative. 
Thus, \eqref{eq:model} captures nearly all the complexity associated with the $x$
coordinate in higher-dimensional kinetic problems with linearized collisions,
and the issue of how best to discretize the speed coordinate is quite independent
of how the other velocity coordinates are discretized.
As an example of how different discretizations may be applied to the different velocity coordinates,
one may refer to the time-independent problems considered in \cite{landreman}, 
in which 2D velocity space is discretized using 
a pseudospectral $x$ discretization
related to the approach considered here, combined with a Legendre modal discretization in pitch angle.
In a similar manner, time-dependent problems in a 2D velocity space could be solved by combining
the $x$ discretization we consider here with a Legendre modal discretization in pitch angle
(or some other pitch angle discretization.) 
The third velocity coordinate, gyro-angle, is often averaged out of kinetic equations in plasma physics
due to rapid particle gyration in a magnetic field,
but this coordinate too could be included using a tensor product approach if desired.

In this article, we compute high-accuracy solutions using a true
  (Galerkin) spectral method to represent the projected dynamics of
  (\ref{eq:model}).  
We discretize
  velocity space only, integrating the resulting ordinary differential
  equations exactly in time.  
In future work, we plan to adapt
  the methods developed here
  to implement exponential time differencing schemes \cite{kassam} and
  implicit-explicit Runge-Kutta methods \cite{carpenter} for
  the time-evolution of the coupled problem.  However, for the model
  problem (\ref{eq:model}), any timestepping scheme will decouple
  into independent eigenmodes that behave as predicted by standard
  linear stability theory; thus, spatial
  discretization is our focus here.

  In a separate paper \cite{vsck3}, pseudo-spectral methods will be
  developed that preserve the self-adjoint structure of the discrete
  evolution operator, mimicking the Galerkin operator as closely as
  possible.  Performance on coarse grids, which is of great practical
  importance in high-dimensional plasma turbulence codes, will be
  addressed in detail there, along with comparisons with other
  methods. In the present article, we also look at performance on
  coarse grids and comparisons with existing discretization schemes
  (Section \ref{sec:low_res}), but the main emphasis of the paper is
  on questions of convergence, discrete approximation of the
  continuous spectrum, and the effect of domain truncation. Except in
  Section~\ref{sec:low_res}, we perform calculations in
  quadruple-precision arithmetic to better illustrate the connection
  between the discrete mode amplitudes in the eigenbasis of the
  projected evolution operator and the continuous
  spectral transform of the solution \cite{vsck1}, to demonstrate
  that the mode amplitudes in the orthogonal polynomial basis continue
  to decay exponentially to arbitrarily small scales, and to more
  closely monitor the effect of domain truncation.  The results of
  Figures~\ref{fig:evol2}--\ref{fig:err:compare}
  and~\ref{fig:trunc1}--\ref{fig:trunc3} below are similar in
  double-precision, with relative errors only a few times larger than
  machine precision, which is $2^{-52}$ in double precision
  versus $2^{-106}$ in the calculations shown in the figures.

  While errors introduced by physical approximations in a model are
  typically large compared to the level of accuracy we consider here,
  it is important when assessing the validity of the model to be
  confident that the results of a numerical simulation accurately
  represent the equations one has discretized.  An additional
  motivation for high accuracy is that small numerical errors can
  impede computation of damped eigenmodes, which are sometimes
  important for physical understanding \cite{hatch}.  In applications
  with more dimensions, such as gyrokinetic simulations
  \cite{candy1,barnes1}, lower resolution would be sufficient for
  routine simulations.  Typically tens of modes or fewer are used in
  the $x$ coordinate, each with tens of degrees of freedom in the
  pitch-angle coordinate of velocity space, and double-precision
  arithmetic is sufficient.

  We observe a remarkable feature of the non-classical polynomials of
  \cite{landreman,shizgal}: for certain initial conditions, the exact
  solution (computed using an expensive spectral transform approach)
  leaves the subspace but returns arbitrarily closely (i.e.~without
  losing spectral accuracy) to the solution of the projected dynamics.
  The situation is different for classical Hermite polynomials. While
  the exact solution still arrives in the subspace with spectral
  accuracy after some time, the projected dynamics evolves to a
  different location in this time, with an error that decays only
  algebraically as a function of the dimension of the subspace.

  The structure of the article is as follows.  In
  Sections~\ref{sec:orthog} and~\ref{sec:proj}, we review the
  construction of orthogonal polynomials and present our general
  method for computing solutions of \eqref{eq:model} by projecting the
  PDE onto finite dimensional subspaces.  Section~\ref{sec:num} forms
  the core of the paper.  In \S\ref{sec:evol:modes}, we solve Equation
  \eqref{eq:model} for two different initial conditions,
  $U(x,0)=xe^{-x^{2}}$ and $U(x,0)=x^{2}e^{-x^{2}}$. In both cases we
  find that the non-standard polynomials are effective at representing
  the solution, and that they are much more accurate than classical
  polynomials defined by orthogonality conditions on the interval
  $(-\infty,\infty)$.  The improved behavior is explained by comparing
  the mode amplitudes in the eigenbasis of the projected evolution
  operator to the spectral transform of the solution, computed
  as described in \cite{vsck1}. This transform method is used in
  \S\ref{sec:valid}-\S\ref{sec:low_res} as an independent means of validating the accuracy
  of the Galerkin approach.  In \S\ref{sec:eigen}, we compare the
  discrete eigenfunctions of the projected operator with the
  non-normalizable eigenfunctions of the PDE to see in what sense the
  continuous spectrum is being approximated by a discrete one.
  In \S\ref{sec:low_res} we study the performance of the new family
  polynomials for the low grid resolutions commonly used in
  five-dimensional kinetic simulations of magnetized plasmas, and
  compare their accuracy to that of classical Hermite polynomials and
  the finite difference scheme used by popular gyrokinetic solvers
  \cite{kotsch,numata}. We then investigate the benefits and drawbacks
  of truncating the domain $[0,\infty)$ to a finite interval in
    \S\ref{sec:trunc}. We summarize our results and discuss future
    work in \S\ref{sec:conclusion}.

\section{Orthogonal Polynomials}
\label{sec:orthog}

We consider two classes of orthogonal polynomials on the positive
half-line.  The first, discussed by Shizgal \cite{shizgal}
and recently applied in the context of plasma physics simulations by
Landreman and Ernst \cite{landreman}, are orthogonal with respect to
the weight function
\begin{equation}\label{eq:rho:plus}
  \rho(x) = x^\nu e^{-x^2}, \qquad\quad (x>0).
\end{equation}
We find in \ref{sec:float} that $\nu=2$ is the best choice overall in
floating point arithmetic, although part of the calculation is more
accurately performed with $\nu=0$.  Roundoff error aside, $\nu=2$ is
also most natural since the projected dynamics involves an implicit
change of basis to this case. Indeed, as explained in
Section~\ref{sec:proj}, $\rho(x)$ then agrees with the weight function
of the Sturm-Liouville problem associated with the evolution operator
on the right-hand side of (\ref{eq:model}).  Thus, except in
\ref{sec:float}, computations will be performed using the $\nu=2$
polynomials. For the remainder of the paper, this class of
polynomials will be referred to as ``full polynomials.''

The second class of polynomials we consider are orthogonal with
respect to
\begin{equation}\label{eq:rho}
  \rho(x) = x^\nu e^{-x^2}, \qquad\quad (x\in\mathbb{R})
\end{equation}
over the whole line. However, the odd polynomials will then be
discarded and the even ones restricted to $\mathbb{R}_+$, giving a
different basis for $L^2(\mathbb{R}_+;\rho\,dx)$.  This basis
  remains complete since any function on the half-line can be extended
  to the whole line by even reflection, which leads to an expansion in
  even polynomials only.
  With $\nu=0$ we obtain the even Hermite polynomials (scaled to be
  monic), whereas with $\nu=2$ we obtain the odd Hermite polynomials
  divided by $2^{2j+1}x$:
\begin{equation}\label{eq:hermite}
  p_j(x) = \left\{\begin{array}{ll}
      H_{2j}(x)/2^{2j}, & \nu=0 \\
      H_{2j+1}(x)/(2^{2j+1}x), & \nu=2
    \end{array}\right\}, \qquad j=0,1,2,\dots
\end{equation}
In either case, the same $n$-dimensional subspace and projected
dynamics will result when the polynomials are truncated to degree
$2(n-1)$, although floating-point issues make the $\nu=2$ family more
desirable to work with.  For the remainder of the article, polynomials
in this second class will be referred to as ``even polynomials'' or
``classical polynomials.''  We note that it is preferable to discard
the odd polynomials rather than the even ones since the exact solution
$u(x,t)=U(x,t)e^{x^2}$ satisfies $\partial u/\partial x=0$ at $x=0$
for $t>0$ rather than $u(x,t)=0$; this follows from the representation
(\ref{eq:uv}) below and the results of \cite{vsck1} on the behavior of
bounded solutions of $Lu=\lambda u$ near $x=0$, where $L$ is defined
in (\ref{eq:PDE}) below.  One or the other must be discarded as the
odd polynomials are not orthogonal to the even ones on the half-line.

Monic orthogonal polynomials with respect to an arbitrary weight
function $\rho(x)$ may be constructed via the Stieltjes procedure
\cite{gautschi:68,gautschi:82}.  A number of technical challenges
arise, partly due to the rapid decay of $\rho(x)=x^\nu e^{-x^2}$,
which leads to overflow and underflow problems in floating point
arithmetic, and partly due to poor conditioning of the recurrence
\begin{equation}\label{eq:3:term:brief}
  p_{j+1}(x)=(x-a_j) p_j(x) - b_j p_{j-1}(x)
\end{equation}
for small values of $x$, which amplifies roundoff errors.
We overcome overflow and underflow issues by carrying an extra integer
to extend the exponent range of floating point numbers in the
polynomial evaluation subroutine (see \cite{vsck3}), and avoid
ill-conditioning in the recurrence (\ref{eq:3:term:brief}) by
representing the polynomials in product form
\begin{equation}\label{eq:prod}
   p_j(x) =  \prod_{k=1}^j(x-x^{(j)}_k).
\end{equation}
The roots $x^{(j)}_k$ of $p_j(x)$ are the eigenvalues of the symmetric
tridiagonal Jacobi matrix $A_j$ with entries
\begin{equation*}
  (A_j)_{ii} = a_{i-1}, \quad (1\le i\le j), \qquad
  (A_j)_{i,i+1} = (A_j)_{i+1,i} = \sqrt{b_i}, \quad (1\le i<j).
\end{equation*}
The coefficients $a_j$ and $b_j$ as well as the squared
norms $c_j=\|p_j\|^2$ are computed via
\begin{equation*}
  p_0 = 1, \qquad
  c_0 = \la p_0,p_0 \ra, \qquad
  a_0 = \la xp_0,p_0 \ra/c_0
\end{equation*}
and the following recursion for $j=1,\dots,n$:
\begin{equation}\label{eq:recur}
    p_j(x) = (x-a_{j-1})p_{j-1}(x) - b_{j-1}p_{j-2}(x), \quad\;
  \begin{aligned}
    c_j &= \la p_j,p_j\ra, \\
    a_j &= \la xp_j,p_j \ra/c_j,
  \end{aligned} \quad\;
  b_j = \frac{c_j}{c_{j-1}},
\end{equation}
where $b_0p_{-1}(x)$ is taken to be zero when $j=1$.  In the floating
point (as opposed to symbolic) algorithm, the formula for $p_j(x)$ is
replaced by (\ref{eq:prod}),
and the inner products $\la f,g\ra=\int_0^\infty
f(x)\overline{g(x)}\rho(x)\,dx$ are computed using composite Gaussian
quadrature on subintervals $[x_k^\e j,x_{k+1}^\e j]$ with endpoints
taken to be the zeros of $p_j(x)$ together with $x_0^\e j=0$ and
several additional points $x_{j+1}^\e j$, \dots, $x_{j+r}^\e j$ chosen
to integrate the tails of the integrands accurately.  See \cite{vsck3}
for further details about this approach,
and \cite{gautschi:68,gautschi:82} for additional perspective.

\begin{figure}
\begin{center}
\includegraphics[width=\linewidth,trim=0 20 0 20]{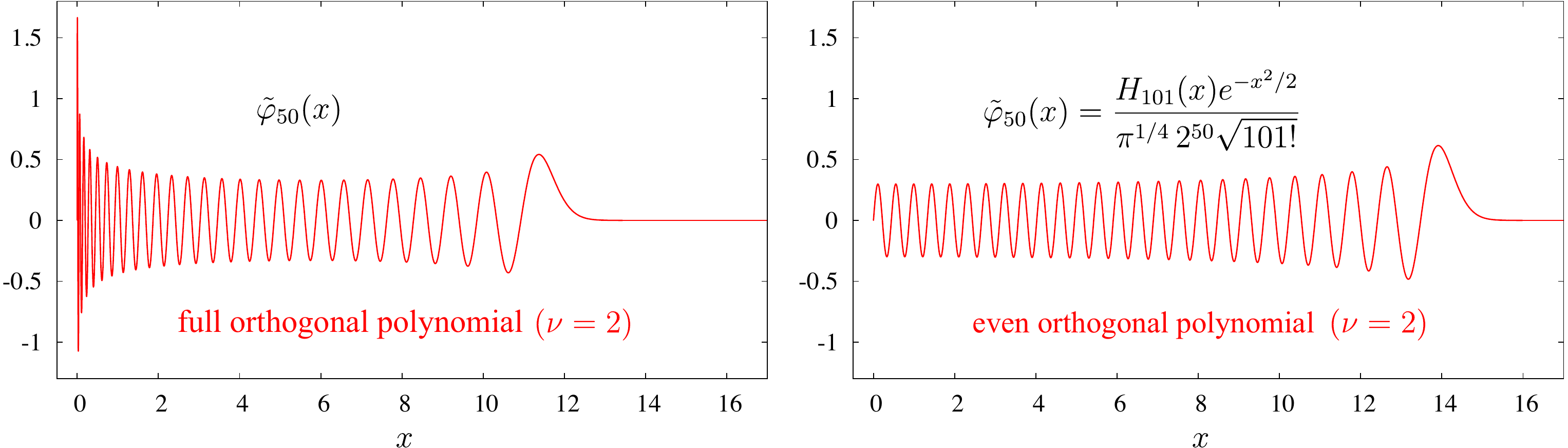}
\end{center}
\caption{\label{fig:basis} Re-scaled basis functions
  $\tilde\varphi_j(x)$ in (\ref{eq:tilde:phi}) with $j=50$ and weight
  function $\rho(x)=x^2e^{-x^2}$.  The main difference between the
  full and even polynomials is that the nodes cluster more tightly at
  $x=0$ in the former case, which causes the oscillations in
  $\tilde\varphi_j(x)$ to be less uniform.  This occurs because $x=0$
  is a true endpoint of the weight function in the full case and is
  merely a symmetry axis in the even case. Note that
  $\tilde\varphi_j(x)$ is odd in the even case due to the factor of
  $x$ in $\rho(x)^{1/2}$.  }
\end{figure}

In what follows, it is useful to define
\begin{equation}\label{eq:varphi}
  \varphi_j(x) = c_j^{-1/2}p_j(x),
\end{equation}
which are (non-monic) polynomials of unit length in
$L^2(\mathbb{R}_+;\rho\,dx)$, and
\begin{equation}\label{eq:tilde:phi}
  \tilde\varphi_j(x) := \rho(x)^{1/2}c_j^{-1/2}p_j(x),
\end{equation}
which are unit vectors in $L^2(\mathbb{R}_+;\,dx)$ that oscillate with
a fairly uniform amplitude between 0 and $x^{(j)}_j$ and then decay
rapidly to zero, as illustrated in Figure~\ref{fig:basis}.
To avoid excessive notation, we use the same symbols $a_j$, $b_j$,
$c_j$, $p_j(x)$, $\varphi_j(x)$ and $\tilde\varphi_j(x)$ in the full
and even cases, with the recurrence relation (\ref{eq:3:term:brief})
replaced, in the even case, by (\ref{eq:even:recur}) below.

\section{Projected Dynamics}
\label{sec:proj}

The first step to solving the PDE (\ref{eq:model}) is to transform it 
to a self-adjoint system.  This is done by defining
\begin{equation}
  u(x,t) = U(x,t)e^{x^2},
\end{equation}
which satisfies
\begin{equation}\label{eq:PDE}
  u_t = -Lu, \qquad
  Lu = -\frac{(\Psi w u')'}{w}, \qquad
  w(x) = x^2e^{-x^2},
\end{equation}
where $u_t:=\partial u/\partial t$ and $u':=\partial u/\partial x$.
Note that if $u$ and $v$ are bounded, $C^2$ functions on $(0,\infty)$,
we have
\begin{equation}
  \langle Lu,v\rangle = \langle u,Lv\rangle, \qquad
  \langle u,v\rangle = \int_0^\infty u(x)\overline{v(x)}\,w(x)\,dx.
\end{equation}
Thus, $L$ is symmetric and densely defined on the Hilbert space
\begin{equation}\label{eq:H}
  \H = L^2(\mathbb{R}_+;w\,dx) =
  \big\{u\;:\;\int_0^\infty |u(x)|^2w(x)\,dx < \infty\big\}.
\end{equation}
Physically, the $k$th moment of the distribution function $U$ may be
computed as the $\H$-inner product of $u$ with $x^k$, or the
$L^2$-inner product of $xe^{-x^2/2}u$ with $xe^{-x^2/2}x^k$:
\begin{equation}\label{eq:moments}
    \int_0^\infty x^2 U(x,t)x^k\,dx = \la u,x^k\ra =
    \int_0^\infty \big(xe^{-x^2/2}u(x,t)\big)\big( xe^{-x^2/2}x^k\big)\,dx.
\end{equation}
It is shown in \cite{vsck1} that $L$ is a singular Sturm-Liouville
operator \cite{coddington, stakgold, krall} on $(0,\infty)$ of limit
circle type at $x=0$ and limit point type at $x=\infty$. The limit
circle case requires a boundary condition, but it suffices to require
that solutions (of $Lu=\lambda u$) remain bounded at $x=0$.  The point
spectrum of $L$ consists of $\lambda=0$ with eigenfunction $u\equiv1$,
and the continuous spectrum is $(0,\infty)$.  A spectral transform
algorithm is developed in \cite{vsck1} that diagonalizes the evolution
operator and expresses the solution as a discrete and continuous
superposition of normalizable and non-normalizable eigenfunctions of
$L$. We will use this computationally expensive method to assess the
accuracy of the projected dynamics below.

In this work, we approximate solutions of the PDE (\ref{eq:PDE})
by projecting onto finite dimensional subspaces of orthogonal
polynomials.  However, it is useful to derive the discrete
  evolution equations without assuming that the basis functions are
  polynomials (or even orthogonal).  Let
$\varphi_0,\dots,\varphi_{n-1}\in \H$ be linearly independent and
consider the subspace $\V=\opn{span}_{0\le j< n}\varphi_j$.  Define
$\Phi:\mathbb{C}^n\rightarrow \V\subset \H$ and its adjoint
$\Phi^*:\H\rightarrow \mathbb{C}^n$ by
\begin{equation}
  (\Phi\vec\alpha)(x) = \sum_{j=0}^{n-1} \alpha_j\varphi_j(x), \qquad
  (\Phi^*u)_i = \langle u,\varphi_i\rangle.
\end{equation}
Let $P=\Phi(\Phi^*\Phi)^{-1}\Phi^*$ be the orthogonal projection from
$\H$ onto $\V$.  The projected dynamics of (\ref{eq:PDE}) onto $\V$ is
given by
\begin{equation}\label{eq:u:tilde}
  \partial_t u_p = -PLu_p,
\end{equation}
where $u_p$ remains in $\V$ for all time.  In weak form, we have
\begin{equation}
  \langle \partial_t u_p, v \rangle = -\langle L u_p, v \rangle,
  \qquad (v\in \V).
\end{equation}
Writing $u_p(x,t) = \sum_{j=0}^{n-1} \alpha_j(t)\varphi_j(x)$, we find that
\begin{equation}\label{eq:Mct}
  M \vec\alpha_t = -K\vec\alpha,
\end{equation}
where
\begin{equation}\label{eq:MK}
  M_{ij} = \langle \varphi_j,\varphi_i\rangle, \qquad
  K_{ij} = \langle L\varphi_j,\varphi_i\rangle =
  \int_0^\infty\Psi(x)\varphi_j'(x)\overline{\varphi_i'(x)}\,w(x)\,dx
\end{equation}
are the mass and stiffness matrices associated with the basis
functions $\varphi_j$, respectively.

The finite-dimensional system of ODEs
  (\ref{eq:Mct})--(\ref{eq:MK}) can be solved by a wide variety of
  time-advance methods.  We solve them {\itshape exactly} here in
  order to focus on the effectiveness of the $x$-discretization,
  without further complications arising from the discretization of
  time.  Due to the linear self-adjoint structure of the equations,
  any other scheme will evolve the eigenmodes independently, in
  accordance with standard linear stability theory.  Its behavior can
  therefore be predicted by replacing $e^{-\lambda t}$ in the spectral
  plots below by $R(-h\lambda)^n$, where $R(z)$ is the stability
  function of the scheme \cite{hairer}, $h$ is the timestep, and $n$
  is the number of steps.

Since the $\varphi_j$ are
linearly independent, $M$ is positive definite and has both a Cholesky
factorization and a square root.  $K$ does as well since it is
positive definite on
\begin{equation}
  \{\varphi_0\}^\perp= \{u\in \H\;:\; \int_0^\infty u(x)w(x)\,dx=0\},
\end{equation}
where we assume $\varphi_0(x)=\text{const}.$ The Cholesky approach is
more convenient for our purposes, so let us write
\begin{equation}\label{eq:R12}
  M = R_1^T R_1, \qquad K = R_2^TR_2, \qquad R := R_2 R_1^{-1} = USV^T,
\end{equation}
where the $R_j$ are upper-triangular and $R_j^T$ is the transpose of
$R_j$ (or the Hermitian transpose if the basis functions are
complex-valued).  Assuming $\varphi_0$ is constant, which is
convenient in practice, the first row and column of $K$ are zero.  The
singular value decomposition $R=USV^T$ solves the eigenvalue problem
for $M^{-1}K$:
\begin{equation}\label{eq:exp:MK}
  (M^{-1}K)(R_1^{-1}V) = (R_1^{-1}V)S^2, \qquad
  e^{-M^{-1}Kt} = (R_1^{-1}V) e^{-S^2t}(V^TR_1).
\end{equation}
The solution of (\ref{eq:u:tilde}) is then
\begin{equation}\label{eq:proj:soln1}
  u_p(x,t)=\Phi(x)\vec\alpha(t), \qquad
  \vec\alpha(t)=e^{-M^{-1}Kt}\vec\alpha(0),
\end{equation}
where $\Phi(x)=(\varphi_0(x),\dots,\varphi_{n-1}(x))$ is treated as a row
vector. Since $PL\Phi = \Phi M^{-1}K$, the eigenfunctions of $PL$ are
the columns of $\Phi(x)R_1^{-1}V$, which are orthonormal in $\H$.

We note that since the constant function $\varphi_0$ is in the basis set,
$P\varphi_0=\varphi_0$ and
\begin{equation*}
  \der{}{t}\int_0^\infty x^2U_p(x,t)\,dx = \partial_t\la u_p,1\ra =
  -\la PLu_p,1\ra = -\la u_p,LP1\ra = 0,
\end{equation*}
i.e.~mass is conserved exactly by the projected dynamics.
The same is true when the domain is truncated
in \S\ref{sec:trunc} since we impose Neumann boundary conditions at
the right endpoint. The constant function remains an eigenfunction
of $L$ in that case.

In floating point arithmetic, computing the SVD of $R$ is more
accurate than forming $K$ and computing its eigenvalues since we avoid
squaring the condition number.  $R_1$ and $R_2$ can be obtained
directly (without forming $M$ and $K$) as follows.  First, we choose a
quadrature scheme $x_k,\mu_k$ such that the matrix entries $M_{ij}$,
$K_{ij}$ in (\ref{eq:MK}) are accurately approximated by
\begin{equation}\label{eq:MK2}
  M_{ij} = \sum_{k=1}^N \varphi_j(x_k)\overline{\varphi_i(x_k)}w(x_k) \mu_k, \quad
  K_{ij} = \sum_{k=1}^N \Psi(x_k)\varphi'_j(x_k)\overline{\varphi'_i(x_k)}w(x_k) \mu_k.
\end{equation}
We choose the $x_k$ and $\mu_k$ corresponding to the composite
quadrature rule used to compute $\la p_j,p_j\ra$ and $\la xp_j,p_j\ra$
in (\ref{eq:recur}) with $j=n$; see \cite{vsck3} for further details.
Note that (\ref{eq:MK2}) may be written
\begin{equation}\label{eq:E12}
  M = E_1^TE_1, \qquad K = E_2^TE_2, \qquad
  \begin{aligned}
    E_{1,kj} &= \sqrt{w(x_k)\mu_k}\,\varphi_j(x_k), \\
    E_{2,kj} &= \sqrt{\Psi(x_k)w(x_k)\mu_k}\,\varphi'_j(x_k).
  \end{aligned}
\end{equation}
$R_1$ and $R_2$ are then obtained by QR-factorization: $E_1=Q_1R_1$,
$E_2=Q_2R_2$. Since the zeroth column of $E_2$ is zero, we actually
perform the QR factorization of columns 1 through $n-1$ to obtain
$\tilde{R}_2$ in (\ref{eq:R}) below.

To reduce roundoff errors in the numerically computed singular values
of $R$, we compute its pseudo-inverse before computing its SVD.  In
more detail, note that since $\varphi_0=\text{const}$,
\begin{equation}\label{eq:R}
  R = R_2R_1^{-1} =
  \left(\begin{array}{c|c} 0 & 0 \\ \hline \rb{-3}{0} & \rb{-3}{\tilde R_2}\end{array}\right)
  \left(\begin{array}{c|c} * & * \\ \hline \rb{-3}{0} & \rb{-3}{\tilde R_1^{-1}}\end{array}\right)
= \left(\begin{array}{c|c} 0 & 0 \\ \hline \rb{-3}{0} & \rb{-3}{\tilde R_2 \tilde R_1^{-1}}\end{array}\right)
= USV^T,
\end{equation}
where $\tilde R_j$ is obtained from $R_j$ be deleting the zeroth row and column.
It follows that
\begin{equation}\label{eq:pinv}
  \opn{pinv}(R) = \left(\begin{array}{c|c} 0 & 0 \\ \hline \rb{-3}{0} &
      \rb{-3}{\tilde R_1 \tilde R_2^{-1}}\end{array}\right) = V\!\opn{pinv}(S)U^T.
\end{equation}
It is well-known \cite{Achain,demmel,cox:higham} that error bounds on the
numerically-computed SVD of an $n\times n$ matrix $A$ are of the form
$Cn^2\veps\|A\|$, where $C$ is a constant (independent of $A$ and $n$);
$\veps$ is machine precision; and $\|\cdot\|$ is either the 2-norm or
the Frobenius norm, depending on whether Givens rotations or
Householder reflections are used in the bidiagonal reduction
process. Thus, we expect that computing $V$ and $S$ via
(\ref{eq:pinv}) will give better results than via (\ref{eq:R}) if
$\|\!\opn{pinv}(R)\|\ll\|R\|$.
We will see in \ref{sec:float} that this is the case when
$\varphi_0,\dots,\varphi_{n-1}$ are orthogonal polynomials with
respect to the weight $w(x)$.

We also note that since $\tilde R_1$ and $\tilde R_2$ are
upper-triangular, applying $\tilde R_2^{-1}$ to $\tilde R_1$ from the
right causes information to propagate from left to right along the
rows of $\tilde R_1$.  More precisely, the first $j$ columns of
$\tilde R_1\tilde R_2^{-1}$ do not depend on columns $k>j$ of $\tilde
R_1$ or $\tilde R_2$.  Thus, numerical error in high-index basis
functions does not corrupt more accurately computed low-index basis
functions in the initial phase of computing $\opn{pinv}(R)$ before
performing the SVD.  Sources of numerical error that tend to be larger
for high-index basis functions include the process of constructing the
orthogonal polynomials $\varphi_j(x)$, the evaluation of derivatives
in the corresponding columns of $E_2$, and quadrature error in the
formulas $M_{ij}=E_{1,:i}^TE_{1,:j}$ and $K_{ij}=E_{2,:i}^TE_{2,:j}$,
where e.g.~$E_{1,:j}$ denotes the $j$th column of $E_1$.

\section{Numerical Results}
\label{sec:num}

We now study the accuracy of approximating solutions of $u_t=-Lu$ by
their projected dynamics.
To illustrate typical behavior, we study two initial
conditions $u(x,0)=f_j(x)$, namely
\begin{equation}\label{eq:ex12:again}
  \text{Example 1:}\quad f_1(x)=x, \qquad\quad
  \text{Example 2:}\quad f_2(x)=x^2.
\end{equation}
Note that $Lf_1$ has a singularity at the origin, so $u_t(x,0)$ blows
up as $x\rightarrow0$ in that case.  For this reason, Example 2 is
easier and presented first below.  We did not switch the labels since
it is easier to remember that Example $k$ corresponds to $f(x)=x^k$.

The cause of the singularity in Example 1 is that $f(x)$ must
be even in order to represent a smooth, radially symmetric function
in three-dimensional velocity space. 
The reader may wonder about the physical relevance of such
an initial condition. Perhaps
surprisingly, this case is indeed physically relevant and occurs in
several practical situations. For example, it occurs when one computes
the resistivity of a plasma (see Section 3 in reference
\cite{landreman}) by solving the time-dependent kinetic equation for
the distribution function until a steady-state is
reached.

We write the solution (\ref{eq:proj:soln1}) in the form
\begin{equation}\label{eq:up:alpha}
  u_p(x,t) = \Phi(x)\vec\alpha(t), \qquad
  \vec\alpha(t) = Ve^{-S^2t}V^T \vec\alpha,
\end{equation}
where $\vec\alpha$ is short for $\vec\alpha(0)$, $V$ and $S$ were
defined in (\ref{eq:R12}), and the matrix $R_{1}$ does not appear in
the second expression because $\nu=2$. Note however that, as discussed
in \ref{sec:float}, the accuracy is improved if $\Phi(x)$ is corrected
for loss of orthogonality by applying $R_1^{-1}$ from the right, where
$R_1$ was defined in (\ref{eq:R12}) and deviates from the identity due
to orthogonality drift in floating point arithmetic. In the numerical
results that follow, $\varphi_k$ stands for the $k$th column of
$\Phi(x)R_1^{-1}$.  See (\ref{eq:proj:soln2}), (\ref{eq:alpha:beta}),
and the last paragraph of \ref{sec:float} for further clarification of
how (\ref{eq:up:alpha}) is computed.

The vector $\vec\alpha$ can be computed analytically for the two
examples in (\ref{eq:ex12:again}).  In the full (as opposed to even)
case, we have
\begin{equation}
  x = p_1 + a_0, \qquad
  x^2 = p_2 + (a_0+a_1)p_1 + (a_0^2+b_1).
\end{equation}
Setting $p_j=\sqrt{c_j}\varphi_j$ and computing $a_j$, $b_j$,
$c_j$ via (\ref{eq:recur}) yields
\begin{equation}\label{eq:alpha:full}
  \begin{aligned}
  x &= \left(\frac{1}{\pi^{1/4}}\right)\varphi_0 +
  \left(\frac{\sqrt{6\pi-16}}{4\pi^{1/4}}\right)\varphi_1, \\
  x^2 &= \left(\frac{3}{4}\pi^{1/4}\right)\varphi_0 +
  \left(\frac{\pi^{1/4}}{\sqrt{6\pi-16}}\right)\varphi_1
  + \left(\frac{\pi^{1/4}\sqrt{9\pi-28}}{2\sqrt{6\pi-16}}\right)\varphi_2,
  \end{aligned}
\end{equation}
which give the coefficients $f(x)=\sum \alpha_k\varphi_k$.

In the even case, we note that $h_j(x)=H_j(x)2^{-j}$ are monic
orthogonal polynomials with weight function $e^{-x^2}$ on $\mathbb{R}$
and satisfy
\begin{equation*}
  h_0(x) = 1, \qquad h_1(x) = x, \qquad
  h_{j+1}(x) = xh_j(x) - \frac{j}{2}h_{j-1}(x).
\end{equation*}
The functions $p_j(x) = h_{2j+1}(x)/x$ are then orthogonal with respect
to $w(x)=x^2e^{-x^2}$ on the half-line and satisfy
\begin{equation}\label{eq:even:recur}
  p_0(x) = 1, \qquad p_1(x) = x^2-a_0, \qquad
  p_{j+1}(x) = (x^2-a_j)p_j(x) - b_jp_{j-1}(x),
\end{equation}
where $a_j=\big(2j+\frac{3}{2}\big)$ and $b_j =
j\big(j+\frac{1}{2}\big)$. It follows that $c_j=\|p_j\|^2$ is given by
\begin{equation*}
  c_j = \int_0^\infty p_j(x)^2w(x)\,dx =
  c_0\prod_{i=1}^j b_i = \frac{\sqrt{\pi}(2j+1)!}{2^{2j+2}}.
\end{equation*}
From $x=\sum_j\la x,p_j\ra c_j^{-1/2}\varphi_j$ and $x^2=\big(a_0\sqrt{c_0}\varphi_0 +
\sqrt{c_1}\varphi_1\big)$ we obtain
\begin{equation}\label{eq:alpha:even}
  x = \sum_{j=0}^\infty \frac{(-1)^{j+1}(2j-3)!!}{\pi^{1/4}\sqrt{(2j+1)!}} \varphi_j, \qquad
  x^2 = \left(\frac{3}{4}\pi^{1/4}\right)\varphi_0 + \left(\frac{\sqrt{6}\pi^{1/4}}{4}\right)\varphi_1.
\end{equation}
An intermediate step is to show that $\la x,p_j\ra =
-h_{2j}(0)/(4j-2)$.  Note that $(-3)!!=-1$ and $(-1)!!=1$ in
(\ref{eq:alpha:even}).  For large $j$, Stirling's approximation gives
\begin{equation}
  \frac{(2j-3)!!}{\sqrt{(2j+1)!}}\approx
  \frac{j^{-7/4}}{(64\pi)^{1/4}}\left(1 + \frac{3}{16j} +
\frac{105}{512j^2}+\cdots\right),
\end{equation}
so the coefficients $\alpha_j$ decay slowly in the even case of
Example 1.

\begin{figure}[t]
\begin{center}
\includegraphics[width=.9\linewidth,trim = 0 20 0 20]{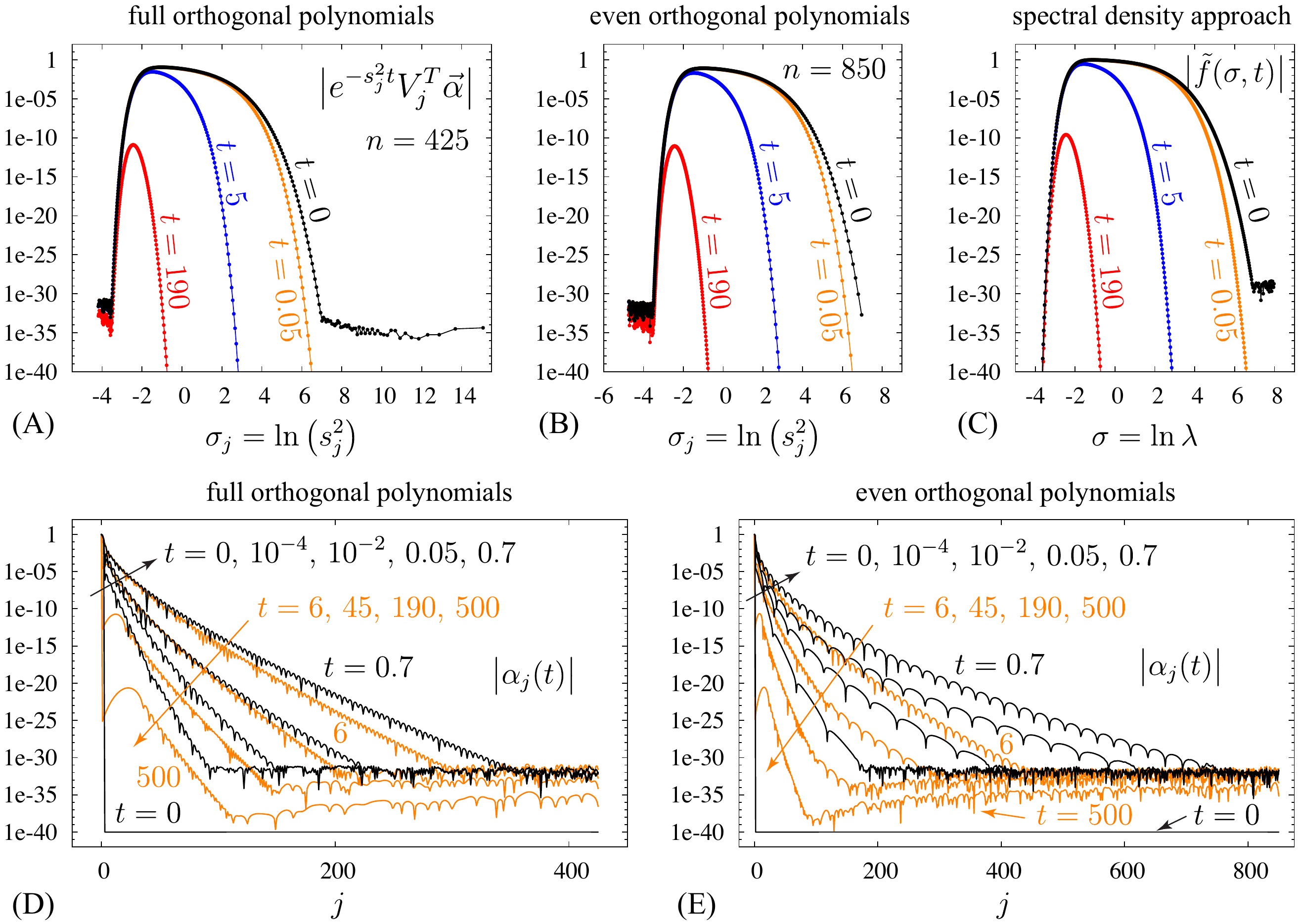}
\end{center}
\caption{\label{fig:evol2} Evolution of mode amplitudes for Example 2
  in the eigenbasis (A,B); the spectral density representation (C);
  and the orthogonal polynomial basis (D,E). While both sets of
  orthogonal polynomials are able to represent the solution of Example
  2, the even polynomials require twice as many basis functions to
  achieve the same level of accuracy. The $s_j=0$ mode is excluded
  in (A,B). }
\end{figure}

\subsection{Evolution of mode amplitudes in the eigenbasis and polynomial basis}
\label{sec:evol:modes}

Since the columns of $\Phi(x)V$ form an orthonormal set of eigenfunctions
for $PL$ on $\V$, the components of the vector
\begin{equation}
  e^{-S^2t}V^T\vec\alpha
\end{equation}
represent the mode amplitudes of $u_p(x,t)$ in the eigenbasis while
the components of $\vec\alpha(t)$ in (\ref{eq:up:alpha}) represent the
mode amplitudes in the orthogonal polynomial basis.  Panels (A,B,D,E)
of Figure~\ref{fig:evol2} show the evolution of both sets of mode
amplitudes for Example 2 in the full and even cases.  For comparison,
we also plot the spectral transform of the solution in (C), which was
computed at 512 equally spaced points between $\sigma=-4$ and
$\sigma=8$ using the algorithm described in \cite{vsck1}.  Here
$\sigma=\ln\lambda$ is the spectral parameter used in \cite{vsck1} to
represent the solution as a discrete and continuous superposition of
eigenfunctions. In more detail, the solution of the PDE in the
infinite dimensional space $\H$ may be written
\begin{equation}\label{eq:uv}
  u(x,t) = \frac{4}{\sqrt{\pi}}\hat{f}(0) + 
  \int_{-\infty}^\infty 
  v(x,\sigma)\tilde{f}(\sigma,t)\,d\sigma,
\end{equation}
where $v(x,\sigma)=u_1(x,e^\sigma)/Y(e^\sigma)$, $u_1(x,\lambda)$ is a
solution of $Lu=\lambda u$ with appropriate boundary conditions at
$x=0$, $Y(\lambda)$ is a scale factor defined in \cite{vsck1},
\begin{equation}
  \tilde f(\sigma,t) = \hat f(e^\sigma)
    e^{-e^\sigma t}\,Y(e^\sigma)\rho'(e^\sigma)e^\sigma, \qquad
    \hat f(\lambda) = \int_0^\infty f(x)u_1(x,\lambda)w(x)\,dx,
\end{equation}
and $\rho'(\lambda)$ is the spectral density function associated with
the singular Sturm-Liouville problem
\cite{coddington,fulton2008a,vsck1}.  In order to understand how the
results in panels A and B are related to those of panel C, consider
the following. Writing $V_j$ for the $j$th column of $V$ and
$S=\opn{diag}(s_j)$, the solution of the projected dynamics is given by
\begin{equation}\label{eq:gamma}
  u_p(x,t) = \sum_j v_j(x)\gamma_j(t), \qquad
  v_j(x) = \Phi(x)V_j, \qquad \gamma_j = e^{-s_j^2t}V_j^T\vec\alpha.
\end{equation}
Comparing (\ref{eq:uv}) and (\ref{eq:gamma}), we can interpret the
projected dynamics solution as having one component that represents
$(4/\sqrt\pi)\hat f(0)$ and the others approximating the integral
$\int v(x,\sigma)\tilde f(\sigma,t)\,d\sigma$.  Since the spacing of
$\sigma_j=\ln(s_j^2)$ is not uniform and the functions $v_j(x)$ are
not normalized in the same way as $v(x,\sigma_j)$, the vertical
scaling in panels (A), (B) and (C) is not expected to be the same.
Nevertheless, these graphs are remarkably similar and give insight
into how the projected dynamics mimic the continuous dynamics.

Panels (D) and (E) of Figure~\ref{fig:evol2} show the evolution of the
coefficients $\alpha_j(t)$ in the orthogonal polynomial basis in the
full and even cases.  We added $10^{-40}$ to all the coefficients to
keep the $t=0$ modes visible in the plot.  In both cases, the number
of active modes grows from 3 or 2 at $t=0$ to 400 or 800 at
$t\approx0.7$, then decays down to one mode (the steady state) as
$t\rightarrow\infty$.  The curves are plotted in black or orange
depending on whether the number of active modes is growing or decaying
at that time.  Panels (D) and (E) show that twice as many modes are
necessary to represent the solution using even polynomials instead of
full polynomials.  Panels (A) and (B) show why this is the case: the
full polynomials are more efficient at sampling the interval
$-4\le\sigma\le 8$ of interest.  Indeed, the eigenvalue distribution
in panel (B) is heavily skewed to over-sample the low end of this
spectral window.

\begin{figure}[t]
\begin{center}
\includegraphics[width=.9\linewidth,trim=0 20 0 20]{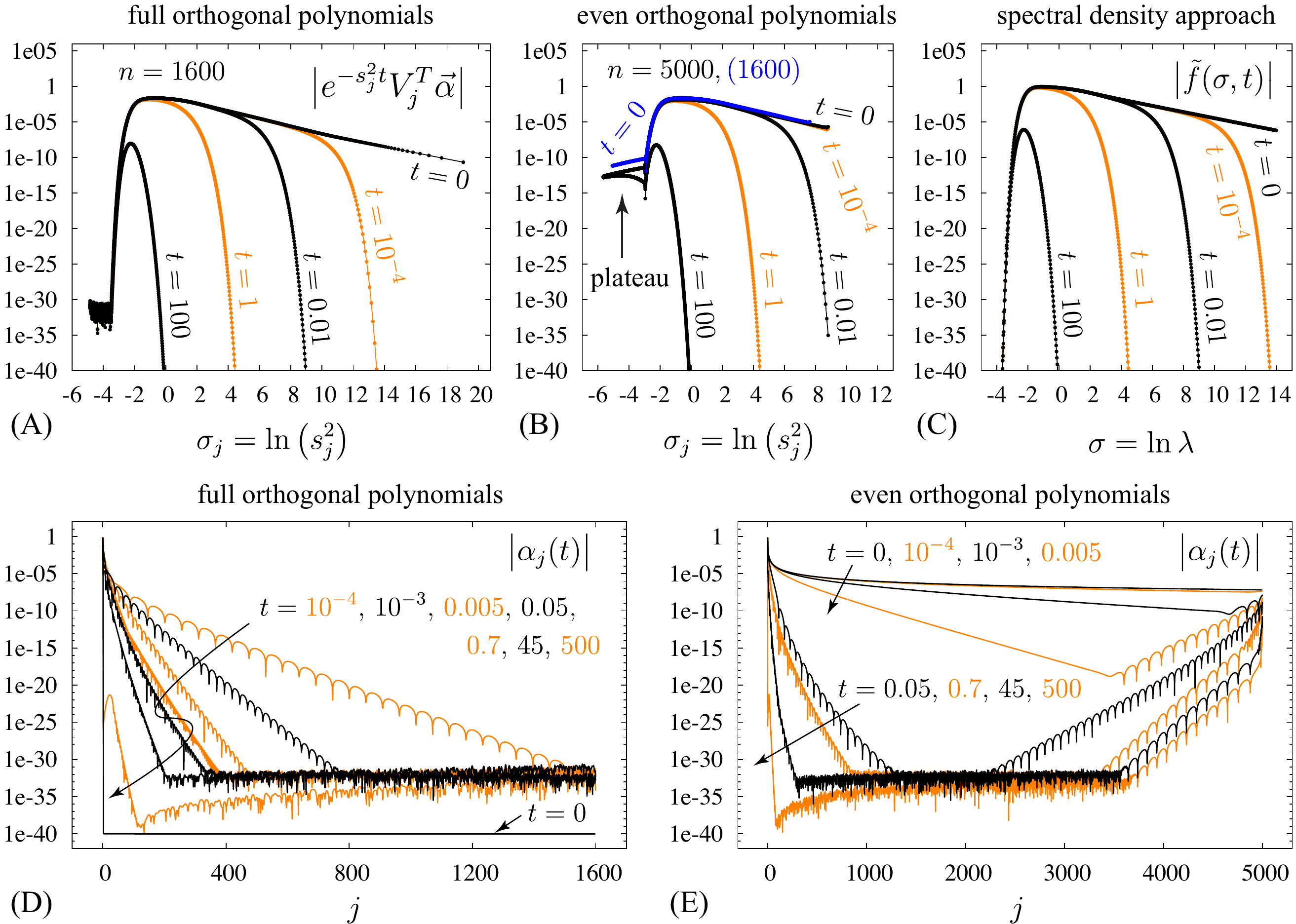}
\end{center}
\caption{\label{fig:evol1} Evolution of mode amplitudes for Example 1
  in the eigenbasis (A,B); the spectral density representation (C);
  and the orthogonal polynomial basis (D,E).  Large errors in the mode
  amplitudes for small values of $s_j$ in panel (B) are not damped out
  quickly by $e^{-s_j^2t}$, and therefore persist for large times in
  panel (E).  In panel (D), the mode amplitudes $|\alpha_j(t)|$
  decrease until $t=0.05$ or so, then enter a brief growth period
  until $t=0.7$, then decay again toward steady-state.
}
\end{figure}

Figure~\ref{fig:evol1} shows that for Example 1 the full polynomials
are again able to mimic the behavior of the continuous problem in the
sense that the mode amplitudes in the eigenfunction basis closely
resemble the spectral transform of the initial condition. Moreover, for
$t\ge10^{-4}$, the mode amplitudes in the orthogonal polynomial basis
decay exponentially to roundoff error using 1600 modes.  By contrast,
in the even case, the mode amplitudes in the eigenfunction basis reach
a plateau when $\sigma$ decreases below -3, and cease to resemble the
spectral transform $\tilde f(\sigma,0)$.  As a consequence, the mode
amplitudes in the orthogonal polynomial basis do not decay to roundoff
error accuracy until $t$ reaches $0.01$ or so, and even then the mode
amplitudes $\alpha_j(t)$ become large again for large $j$. Thus, the
even polynomials are not able to represent the solution of Example 1
effectively even for large values of $t$.

\begin{figure}[t]
\begin{center}
\includegraphics[width=\linewidth,trim=0 20 0 20]{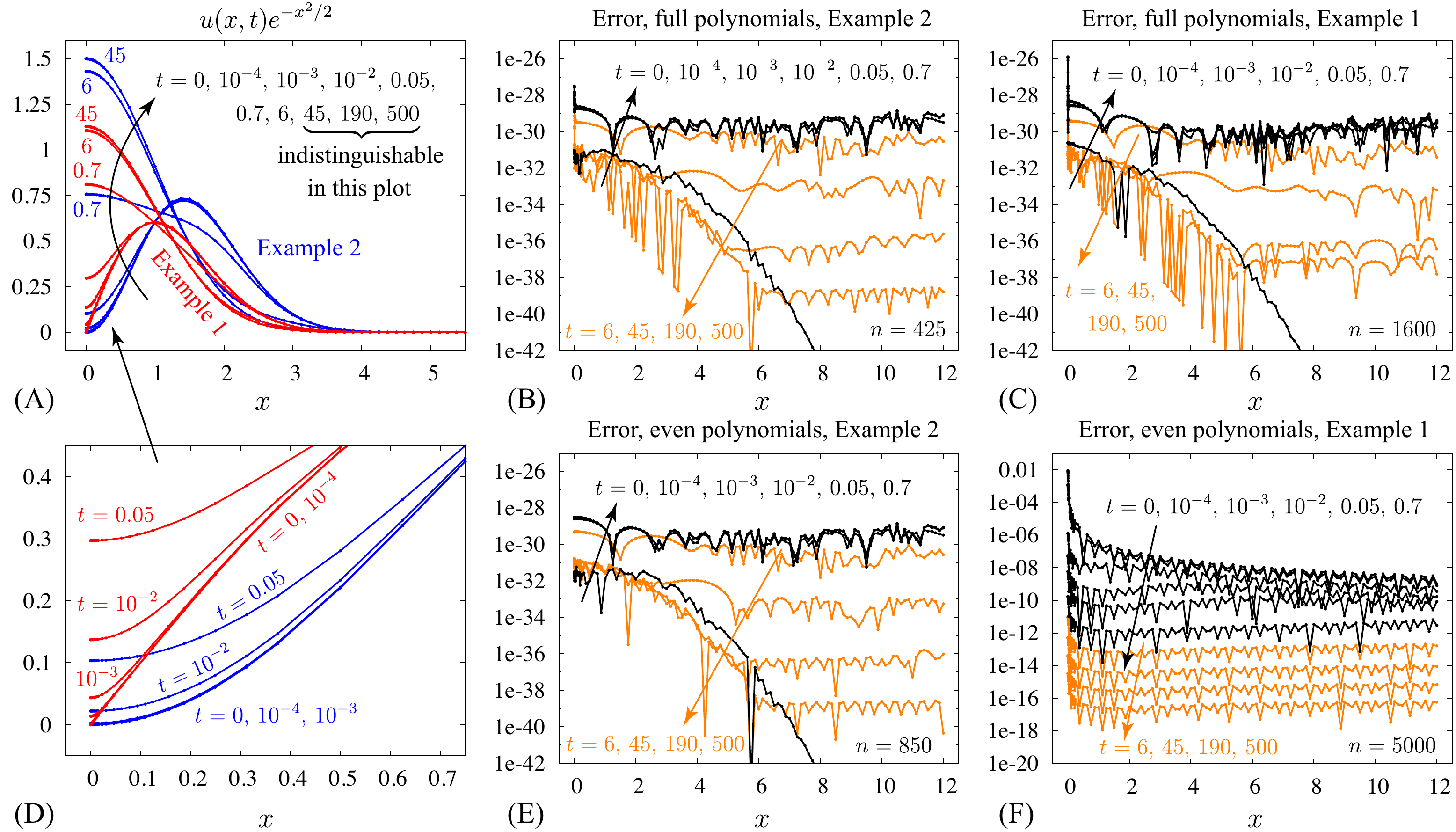}
\end{center}
\caption{\label{fig:err:compare} Evolution and error of projected
  dynamics in spaces of full and even orthogonal polynomials. The
  spectral transform approach of \cite{vsck1} was used for the
  ``exact'' solution.  (B,C,E) both methods reach
  roundoff error in quadruple-precision. (F) the projected
  dynamics with even orthogonal polynomials does not perform well with
  5000 modes until the solution approaches steady-state.}
\end{figure}

\subsection{Validation of accuracy}
\label{sec:valid}

Of course, it is not guaranteed that the full polynomials yield the
correct answer for $t\ge10^{-4}$ just because the mode amplitudes
decay fast enough to reach roundoff error at that point.  The true
solution might leave the subspace $\V=\opn{span}\varphi_k$ initially
and come back to a different point in the subspace than the projected
dynamics predicts. This was the motivation for developing the spectral
transform approach in \cite{vsck1}, where we know the analytic form of
the solution and can use it as an independent check of the correctness
of the projected dynamics. In Figure~\ref{fig:err:compare}, we plot
the solution, scaled by $e^{-x^2/2}$, for Examples~1 and~2 at several times,
together with
the difference between the solution of the projected dynamics and the
one obtained using the spectral transform approach, both scaled by
$e^{-x^2/2}$.  We considered other scalings, namely
$e^{-x^2}$ (since $U=ue^{-x^2}$) and also $xe^{-x^2/2}$; the scaling
$e^{-x^2/2}$ is the one for which the plots are the most readable. We
added $(10^{-32-2x})$ to the error plots to avoid evaluating
$\log_{10}(0)$ when the two methods agree to all 32 digits.

At $t=0$, both methods represent the solution with high relative
accuracy for Example 2, but only the full polynomials are able to do
this for Example 1.  This is because the only errors at $t=0$ in
panels (B,C,E) of Figure~\ref{fig:err:compare} are roundoff errors in
the 2--3 nonzero coefficients in (\ref{eq:alpha:full}) and
(\ref{eq:alpha:even}), while in panel (F) the series in
(\ref{eq:alpha:even}) has been truncated at 5000 terms, leading to
much larger errors at $t=0$.  In panels (B) and (E), which correspond
to Example 2, the two approaches agree to more than 28 digits of
accuracy for all positive times, though twice as many modes are needed
to achieve this accuracy with even polynomials. However, the errors
are now absolute errors rather than relative errors. This loss of
relative accuracy is expected as the spectral density approach leads
to oscillatory integrals involving a large amount of cancellation to
evaluate $u(x,t)e^{-x^2/2}$ for $t>0$ and $x\ge6$.  Similarly, the
orthogonal polynomial approach yields a superposition $\sum_j
\alpha_j(t)\varphi_j(x)e^{-x^2/2}$ that involves a large amount of
cancellation of digits for $t>0$ and $x\ge6$.  For example,
$\varphi_j(x)e^{-x^2/2}$ remains $O(1)$ well past $x=30$ when $j=425$
while $|u(x,t)e^{-x^2/2}|$ is less than $10^{-30}$ when $x=12$ and
$t\ge0$.  For very large $t$, both methods regain high relative
accuracy since the oscillatory part of the calculation becomes
negligible compared to the steady-state zeroth mode. Panel (C)
confirms that the full polynomials yield small (absolute) errors in
Example 1 for $t\ge10^{-4}$ even though the true solution leaves the
subspace $\V$ for $0<t<10^{-4}$ (due to slower decay of mode amplitudes
$|\alpha_j(t)|$ in Figure~\ref{fig:evol1} than those shown for
$t=10^{-4}$).
Panel~(F) of Figure~\ref{fig:err:compare} shows that this is not true for
even polynomials.

There are two issues in play causing the even polynomials to be less
accurate. It takes longer (by a factor of 66)
for the exact solution (computed using the spectral transform method)
to reach the subspace --- we have checked that the distance from the
exact solution to the 5000-mode even polynomial subspace does not drop
below $10^{-30}$ until $t=0.0066$ --- and when it arrives, the
projected dynamics solution has evolved to a different location in the
subspace, with an error of $1.302\times10^{-8}$ in the $\H$-norm at
$t=0.0066$ with 5000 modes. For comparison, the $\H$-norm of the exact
solution at this time is $0.66354$.

\begin{figure}[t]
\begin{center}
\includegraphics[width=\linewidth,trim=0 20 0 20]{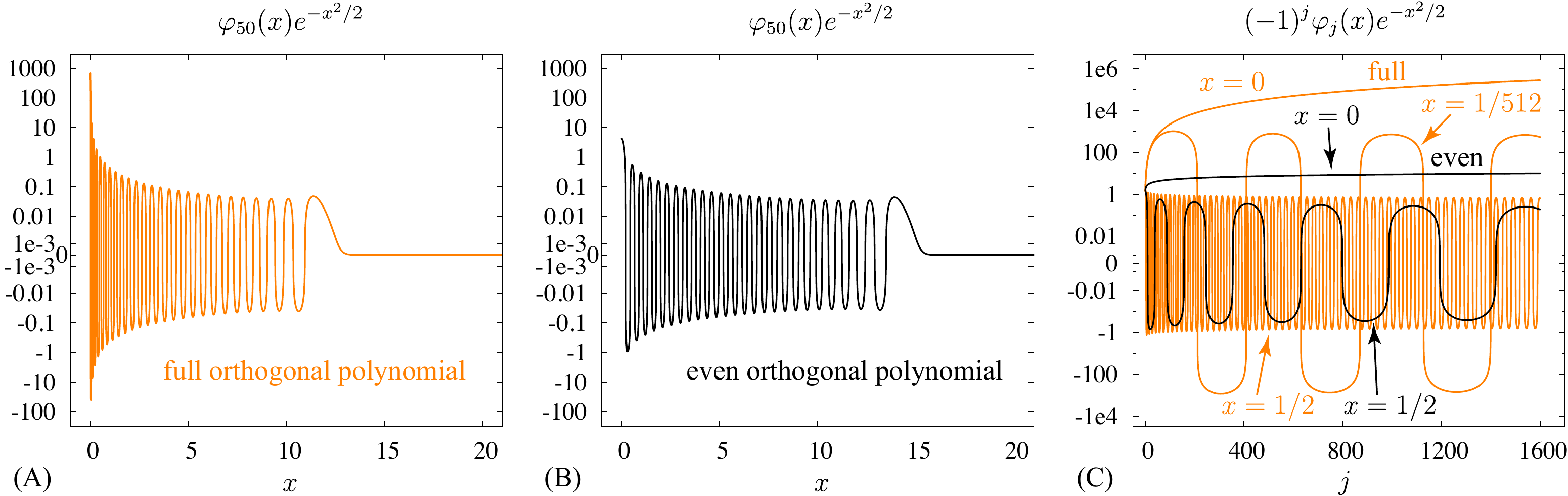}
\end{center}
\caption{\label{fig:growthP0} Full orthogonal polynomials
  $\varphi_j(x)$ (orange) grow much faster with $j$ at $x=0$ than even
  orthogonal polynomials (black).  In panel (C), the factor of
  $(-1)^j$ is included to account for the fact that $\varphi_j(0)$
  alternates in sign. The function $(-1)^j\varphi_j(x)e^{-x^2/2}$
  increases monotonically in $j$ when $x=0$ is held fixed, but becomes
  oscillatory in $j$ for any other fixed value of $x$.  }
\end{figure}

In Figure~\ref{fig:err:compare}(B,C), the errors in the full
polynomial case are much larger near $x=0$ than elsewhere.  This
  is not a cause for concern since physical quantities, such as the
  moments of the distribution function in (\ref{eq:moments}), carry an
  additional factor of $x^2$ in the integrand, which suppresses these
  errors. Nevertheless, it is instructive to identify their source.
  These errors occur because the coefficients $\alpha_j$ in
Figures~\ref{fig:evol2}(D) and~\ref{fig:evol1}(D) carry roundoff
errors that are amplified by the large values of the basis functions
$\varphi_j(x)$ near $x=0$ for large $j$.  As shown in
Figure~\ref{fig:growthP0}, $\varphi_j(0)$ is already close to 1000 in
the full case when $j=50$, and grows to $2.8\times 10^5$ at $j=1600$.
For other values of $x$, $(-1)^j\varphi_j(x)e^{-x^2/2}$ is oscillatory
in $j$, with higher frequency and smaller amplitude oscillations for
larger $x$.  Thus, these functions are only large near $x=0$.  This
growth in the basis functions near $x=0$ occurs because the weight
function $w(x)=x^2e^{-x^2}$ approaches 0 as $x\rightarrow0^+$.
Indeed, we saw in Figure~\ref{fig:basis} that when $j$ is large,
$x\varphi_j(x)e^{-x^2/2}$ oscillates with a fairly uniform amplitude
over a large distance before eventually decaying to zero as
$x\rightarrow\infty$.  In the even case, $x\varphi_j(x)e^{-x^2/2}$
resembles a sine function near $x=0$; hence, $\varphi_j(x)e^{-x^2/2}$
resembles a sinc function.  However, in the full case, the zeros of
$\varphi_j(x)$ are more tightly clustered near $x=0$ since it is a
true endpoint rather than a symmetry axis. The higher
oscillation rate near $x=0$ causes the peaks of
$x\varphi_j(x)e^{-x^2/2}$ to be amplified more when divided by $x$ to
obtain $\varphi_j(x)e^{-x^2/2}$, since $x$ is smaller at the
peaks. Thus, $\varphi_j(0)$ grows more rapidly with $j$ in the full
case than in the even case.


\subsection{Eigenfunctions}
\label{sec:eigen}

Next we consider the connection between eigenfunctions of the
projected operator $PL$ and solutions of $Lu=\lambda u$, which are not
normalizable but serve as basis functions to represent the solution of
$u_t=-Lu$ in a continuous superposition via the spectral transform
(\ref{eq:uv}) and (\ref{eq:gamma}).  Since the exact solution has this
form, one might expect that the accuracy of a discrete approximation
of the continuous spectrum would be limited by the degree to which
these eigenfunctions can be approximated. We find below that this is
not the case: the eigenfunctions need only agree to 2--3 digits for
the time dependent solutions constructed from them to agree to 30 digits.

\begin{figure}[p]
\begin{center}
\includegraphics[width=.95\linewidth,trim=0 20 0 20]{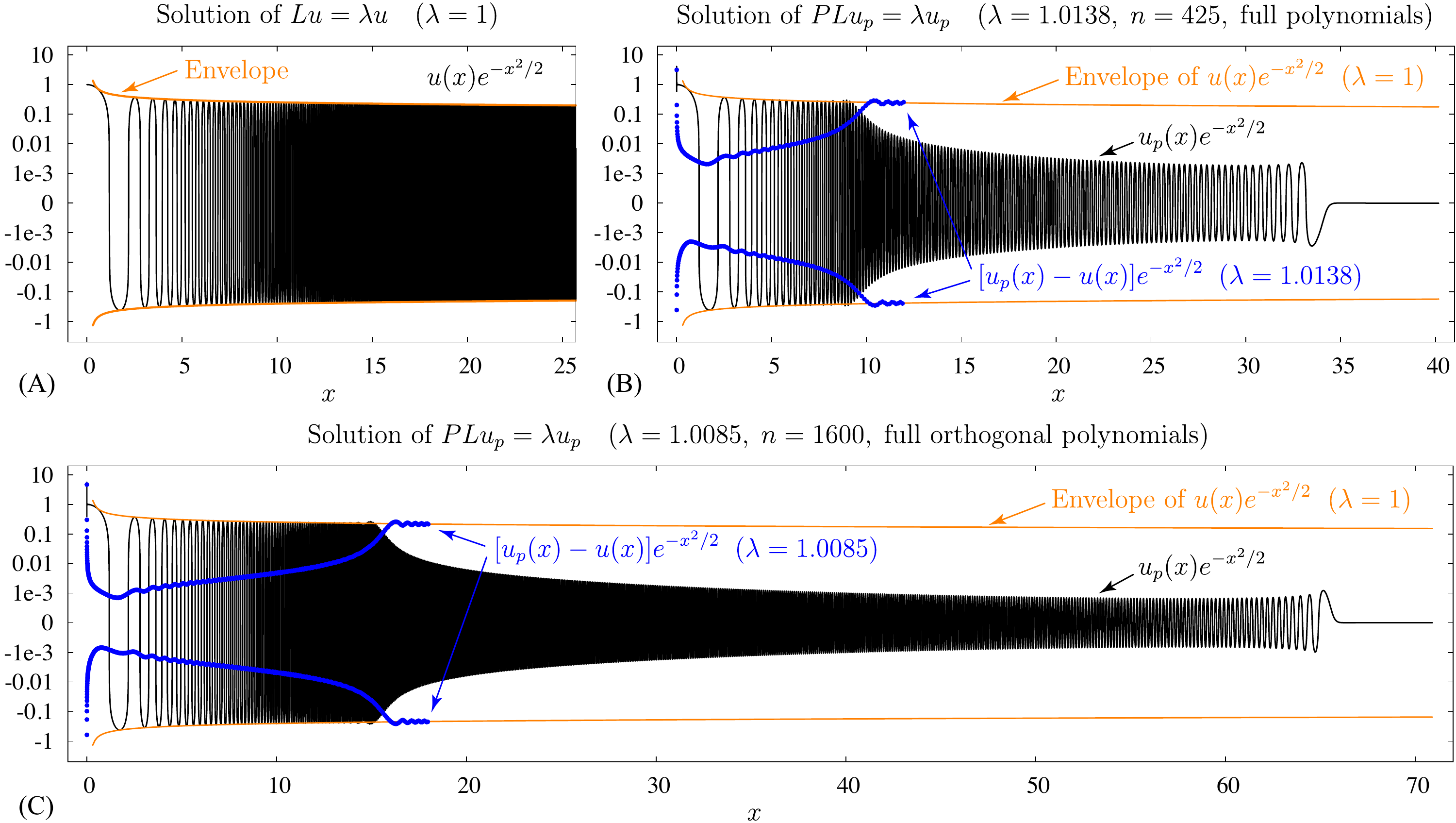}
\end{center}
\caption{\label{fig:efun} Comparison of eigenfunctions of $PL$ and
  solutions of $Lu=\lambda u$.  The asymptotic formula (\ref{eq:ufit})
  was used for the envelopes (orange curves) while a 15th order
  spectral deferred correction scheme was used to solve $Lu=\lambda u$
  to plot the black curve in (A) and the blue markers in (B), (C). }
\end{figure}

Figure~\ref{fig:efun} shows the
solution $u$ of $Lu=\lambda u$ with $u(0)=1$ and $\lambda=1$, as well as
the eigenfunction $u_p$ of $PL$ with eigenvalue closest to 1 for two
choices of $n=\dim\V$.  One boundary condition on $u$ is
sufficient at $x=0$ since the other linearly independent solution
blows up there.  The orange curves show the envelope of the solution
of $Lu=\lambda u$, which we define as the prefactor $A_0p(x)$ in
the asymptotic formula
\begin{equation}\label{eq:ufit}
  u(x)e^{-x^2/2} = A_0 p(x)\cos[q(x)+\theta_0]+O(x^{-15/4}).
\end{equation}
In \cite{vsck1}, the authors show that
\begin{align}
  p(x) &= \jt x^{-1/4}\left[1 +
    \frac{1}{8x\lambda} + \frac{5}{128x^2\lambda^2} +
    \frac{15}{1024x^3\lambda^3}\right], \\
  q(x) &= \jt \sqrt{2\lambda x^5}\left[
    \frac{2}{5} - \frac{1}{6x\lambda} - \frac{1}{16x^2\lambda^2} +
    \frac{1}{64x^3\lambda^3} + \frac{5}{3072x^4\lambda^4} +
    \frac{7-1152\lambda^4}{20480x^5\lambda^5}\right]. \label{eq:q}
\end{align}
The parameters $A_0=0.443935$ and $\theta_0=-1.615039$ were obtained
by fitting the numerical solution of $Lu=\lambda u$, $u(0)=1$,
$\lambda=1$ to the form (\ref{eq:ufit}) for large $x$. Note that the
amplitude of $u(x)$ decays slowly while the frequency grows
rapidly since $p(x)\sim x^{-1/4}$ and $q(x)\sim x^{5/2}$.

The blue markers in Figure~\ref{fig:efun}(B,C) show the extrema
of the (highly oscillatory) residual
\begin{equation*}
  r(x) = [u_p(x)-u(x)]e^{-x^2/2}.
\end{equation*}
We re-scaled the eigenfunction $u_p$ by hand to minimize the amplitude
of the oscillations in $r(x)$.  Due to the logarithmic scale of the
plot, it is important to use the same $\lambda$ for $u$ and $u_p$ when
computing $r$, but setting $\lambda=1$ is sufficient for plotting the
envelope.  For $n=425$ (panel B of the figure), there is an
interval $0.125\le x\le 6$ where $u_p$ agrees with $u$ to 2--3 digits
of accuracy.  For $x>6$, $r(x)$ grows in amplitude since there are not
enough orthogonal polynomial basis functions for $u_p(x)$ to match the
accelerating frequency of oscillation in $u(x)$. As a result, the blue
markers move outward and begin to oscillate about the envelope curve
as $u$ and $u_p$ pass in and out of phase with each other.  At this
point, the projected eigenfunction $u_p$ drops in amplitude by a few
orders of magnitude and enters a plateau phase where it no longer
reaches the envelope curve but still remains significant in
size. Beyond $x=35$, the eigenfunction finally decays rapidly to zero.
The results for $n=1600$ (panel C) are similar, but $|r(x)|<10^{-2}$
over a larger window $0.025<x<12.7$ and is generally 3--10 times
smaller than in the $n=425$ case. The plateau region also grows from
$10<x<35$ when $n=425$ to $16<x<66$ when $n=1600$.

\begin{figure}[p]
\begin{center}
\includegraphics[width=\linewidth,trim=0 20 0 20]{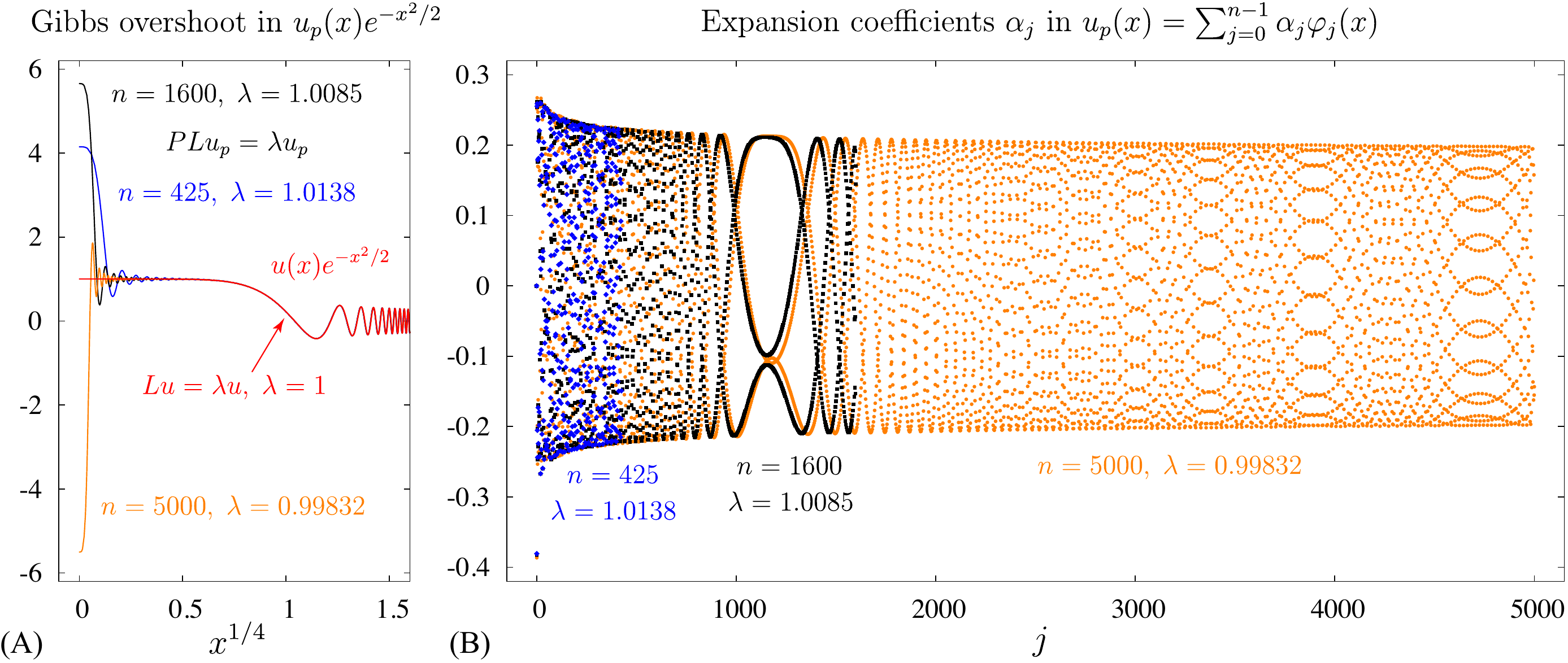}
\end{center}
\caption{\label{fig:gibbsV} Gibbs overshoot phenomenon near $x=0$ (left)
  and coefficients $\alpha_j$ expressing the eigenfunctions $u_p$ in
  the orthogonal polynomial basis. The coefficients do not decay to
  zero as $n$ increases since $\lambda=1$ is part of the continuous
  spectrum of $L$, i.e.~there is no normalizable eigenfunction of $L$
  on $\mathbb{R}_+$ to which $u_p$ could converge.}
\end{figure}

The growth in $r(x)$ as $x\rightarrow0^+$ in
Figure~\ref{fig:efun}(B,C) is due to a Gibbs phenomenon in the
eigenfunctions of $PL$.  This is shown in more detail in
Figure~\ref{fig:gibbsV}(A), which plots the eigenfunction with
eigenvalue closest to 1 for $n=425$ (blue), $n=1600$ (black), and
$n=5000$ (orange), along with the solution $u$ of $Lu=\lambda u$,
$\lambda=1$ (red).  We used a linear scale on the $y$-axis to better
illustrate the magnitude of the overshoot. We also plotted the results
parametrically versus $x^{1/4}$ to better distinguish the oscillations
from the $y$-axis and from each other.  As $n$ increases, the
overshoot becomes taller, narrower, and more oscillatory.

In Figure~\ref{fig:gibbsV}(B), we plot the coefficients when these
eigenfunctions are expanded in the orthogonal polynomial basis.  This
plot illustrates that all of the eigenfunctions $u_p$ of $PL$ are
poorly resolved in the sense that the expansion coefficients do not
decay to zero once $n$ is large enough.  The one exception is the
$\lambda=0$ eigenfunction, $u(x)\equiv1$, which is normalizable and
agrees with $p_0(x)$.  For the others, as $n$ increases, the
higher-frequency polynomials make it possible for $u_p$ to match the
oscillations in $u$ over a greater distance.  Thus, as discussed
above, the window over which $r(x)$ is small grows from $0.125<x< 6$
when $n=425$, to $0.025<x<12.7$ when $n=1600$, to $0.0067<x<21$ when
$n=5000$ (not shown in Figure~\ref{fig:efun}).  The plateau region
also grows with $n$. This is a consequence of all the polynomial basis
functions being present in $u_p$. Indeed, the effective support of
$u_p$ (where $|u_p(x)e^{-x^2/2}|>10^{-30}$, say) is roughly the same
as that of the highest frequency orthogonal polynomial present, which
grows with $n$.  By contrast, we will see below that on a truncated
domain, the eigenfunctions of $PL$ converge to those of $L$, and are
independent of $n$ once $n$ increases beyond the point where the
coefficients $\alpha_j$ have converged to zero.  It is remarkable that
on the infinite domain, in spite of the Gibbs overshoot and rather
poor agreement between eigenfunctions of $PL$ and solutions of
$Lu=\lambda u$, the numerical solution of the projected dynamics
agrees to 30 digits of accuracy to the numerically computed spectral
transform solution, which is built from solutions of $Lu=\lambda u$.

\subsection{Accuracy at low resolution}
\label{sec:low_res}

\begin{figure}
\begin{center}
\includegraphics[width=\linewidth,trim=0 20 0 20]{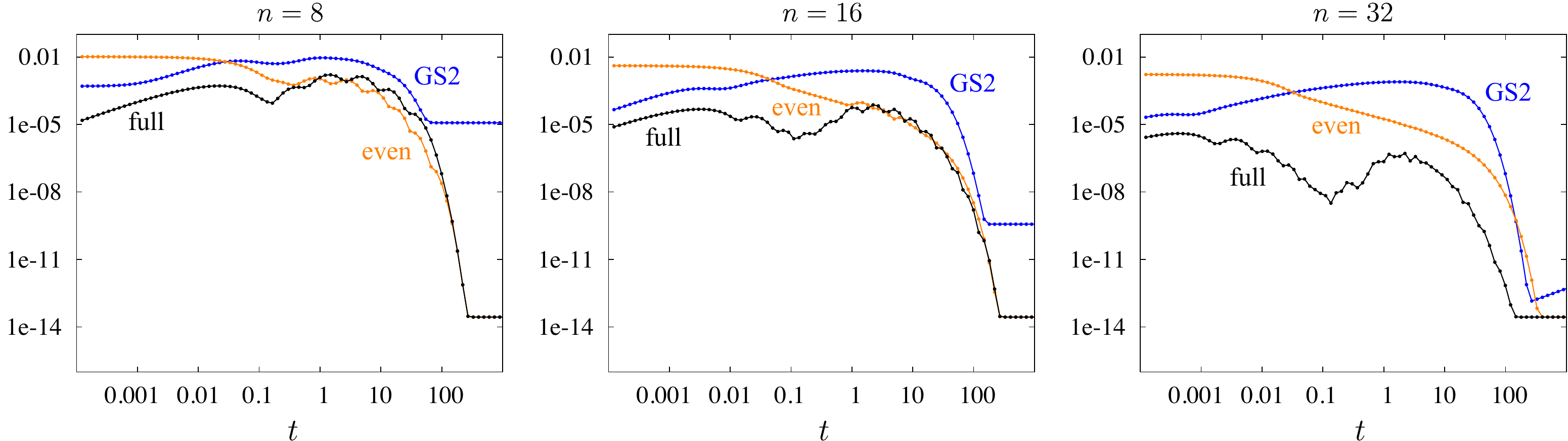}
\end{center}
\caption{\label{fig:smallGrid} Evolution of the error in the Hilbert
  space $\H$ as a function of time for three discretization schemes:
  the full and even polynomials discussed throughout the article, and
  the hybrid Gauss-Legendre/Gauss-Laguerre grid used by the
  gyrokinetic codes GS2 and AstroGK. Here $n$ is the number of modes
  or grid points used to solve (\ref{eq:model}), the initial
  conditions correspond to Example 1, and the spectral transform
  approach of \cite{vsck1} was used for the ``exact'' solution.}
\end{figure}

\begin{figure}
\begin{center}
\includegraphics[width=\linewidth,trim=0 20 0 0]{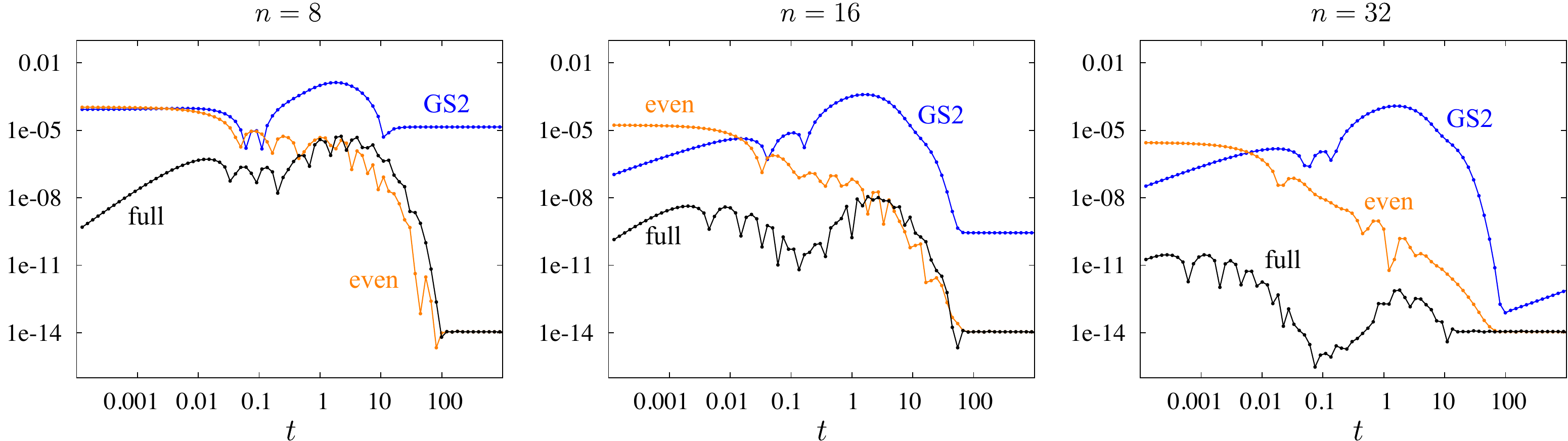}
\end{center}
\caption{\label{fig:smallGridEntropy} Evolution of the error in the
  entropy-like scalar function
  $S(t)=-\int_{0}^{\infty}U^{2}x^{2}e^{x^{2}}$ as a function of time
  for the discretization schemes of Figure~\ref{fig:smallGrid}. As
  before, $n$ is the number of modes or grid points used to solve
  (\ref{eq:model}), the initial conditions correspond to Example
  1, and the spectral transform approach of \cite{vsck1} was used for
  the ``exact'' solution.}
\end{figure}

Numerical simulations of high dimensional kinetic equations are often
so computationally intensive that only a low number of modes or
grid points can be retained to discretize the speed variable. For
example, in the five-dimensional gyrokinetic simulations of plasma
turbulence in magnetic confinement devices, the number of grid points
for the speed variable is typically
8-24 \cite{barnes3}. If the non-classical polynomials studied in this
article are to replace existing discretization schemes for kinetic
simulations, it is important to consider their performance
at low resolution.
In Figures \ref{fig:smallGrid} and \ref{fig:smallGridEntropy}, we
compare the error at low resolution for three different discretization
schemes. The first two are the full and even polynomials considered
throughout the article. The third is that of the popular gyrokinetic
codes GS2 and AstroGK \cite{kotsch,numata}, in which the grid
consists of $n-m$ Gauss-Legendre points on the interval $[0, 2.5]$,
together with $m$ Gauss-Laguerre points on $[2.5,\infty)$, where $m=1$
for $n\leq 12$ and $m=2$ otherwise \cite{barnes3}. In more detail, the
Gauss-Laguerre
points $\tilde x_j\in[0,\infty)$ are transformed to $[2.5,\infty)$ via
$x_j = (6.25+\tilde x_j^2)^{1/2}$, and the weights are scaled by
$w_j=\tilde w_j/(2x_j)$, so that $\int_{2.5}^\infty p(x^2)xe^{-x^2}dx
= \sum_{j=1}^m p(x_j^2)x_je^{-x_j^2}w_j$ whenever $p(\tilde x)$ is a
polynomial of degree less than $2m$.  We actually used $m=\lfloor
n/3\rfloor$ since the errors were smaller.  In GS2 and AstroGK, the
Gauss-Legendre and Gauss-Laguerre points are used for accurate
numerical integration \cite{barnes3}, but derivatives with respect to
the speed variable, such as the ones appearing on the right-hand side
of (\ref{eq:model}), are computed with a 3-point finite-difference
stencil.

Figure \ref{fig:smallGrid} gives the norm of the error in the Hilbert
space $\H$ as a function of time, and Figure
\ref{fig:smallGridEntropy} gives the error in the computed value of
$S(t)=-\int_{0}^{\infty}U^{2}x^{2}e^{x^{2}}dx$ as a function of time.
$S(t)$ is an entropy-like function in the sense that it satisfies
$dS/dt\geq0$, which can be viewed as a ``second principle of
thermodynamics.'' This can be easily seen by multiplying
(\ref{eq:model}) by $x^{2}e^{x^{2}}U$ and integrating from $x=0$ to
$x=\infty$. The initial conditions in Figures~\ref{fig:smallGrid}
and~\ref{fig:smallGridEntropy} are those of the more challenging
Example 1, and the error is measured by comparing the solutions
obtained with the three different discretization schemes with the
solution obtained from the spectral transform approach. These figures
show that for small grid sizes, the full polynomials are several
orders of magnitude more accurate than the even polynomials for small
times, and become comparable at larger times.  For larger grid sizes
($n\ge 32$), the full polynomials remain more efficient than even
polynomials for all time, as seen in the third panel of Figures
\ref{fig:smallGrid} and \ref{fig:smallGridEntropy}, and also in
Figures~\ref{fig:evol2} and~\ref{fig:evol1}. We can also see that at a
fixed grid size, the full polynomials lead to more accurate solutions
than the GS2 scheme for all times, often by several orders of
magnitude.

For the time evolution of the present paper, we evolve to
any time $t$ by diagonalizing the discrete differential operators.
This leads to almost no difference in the running times of GS2 and the
orthogonal polynomial approach. For example, to compute the solution
at all 90 values of $t$ in Figures~\ref{fig:smallGrid}
and~\ref{fig:smallGridEntropy}, the running time (in milliseconds) of
the computational phases of our two codes was
\begin{equation*}
  \begin{array}{c|c|c|c|c|c|c}
    & A\, (n=8)   & B\, (n=8)   &
      A\, (n=16)  & B\, (n=16)  &
      A\, (n=32)  & B\, (n=32)  \\ \hline
      GS2 & 4.64 & 0.076 & 4.30 & 0.121 & 4.59 & 0.233 \\ \hline
      OP &  3.10 & 0.075 & 2.96 & 0.122 & 3.14 & 0.198
  \end{array}
\end{equation*}
The columns labeled $A$ correspond either to balancing the tridiagonal
matrix and computing its eigenvalues and eigenvectors (GS2 approach)
or computing the SVD of $R$ in (\ref{eq:R}), which is precomputed and
read from a file (orthogonal polynomial approach). The same file works
for any $n$ (up to the one used to generate the file) since the $R$
matrices are nested.  The columns labeled $B$ correspond to evolving
the solution from $t=0$ to the times shown in the figures by computing
$Ve^{-\Lambda t}V^Tu_0$, where $V$ is the eigenvector matrix. For
these small values of $n$, the running times depend more on the BLAS
implementation than the number of flops involved, and the first phase
(column A) is actually faster when $n=16$ than $n=8$ in both
approaches. We used the Intel Math Kernel library in a C++ framework
for both codes.

One advantage of the GS2 scheme is that the discretized differential
operator $L$ on the right-hand side of (\ref{eq:PDE}) is tridiagonal.
Thus, if the timestepping scheme of the high-dimensional problem
involves solving one-dimensional subproblems (e.g.~through operator
splitting), the computational cost of applying $L$ or inverting $I+hL$
will be lower for GS2 than for either set of orthogonal polynomials.
Thus, one could potentially use a larger $n$ in GS2 for the same
computational cost as the orthogonal polynomial approach.  However,
increasing $n$ also increases the required storage space and
communication costs, which may be more limiting resources than the CPU
cycles available to solve the local one-dimensional subproblems.
Nevertheless, sparse discretizations are clearly desirable.  We are
developing a banded version of the orthogonal polynomial approach in a
pseudo-spectral framework that is more accurate than GS2 while
retaining much of the sparsity advantage over the dense Galerkin
approach presented here \cite{vsck3}.

\section{Truncation of the Domain}
\label{sec:trunc}

There are three benefits to truncating the domain to a finite
interval. First, the method is easier to implement as the
integration domain in the inner products is fixed. We continue
to use a composite Gaussian quadrature rule using the zeros
of $p_n(x)$ as the endpoints of the integration sub-intervals.
However, it is no longer necessary to deal with the last
sub-interval as a special case since it no longer extends
to infinity.  Second, the coefficients $c_j=\la p_j,p_j\ra$
and $b_j = c_j/c_{j-1}$ grow less rapidly on a truncated domain.
For example, numerical experiments suggest that for large $k$
\begin{equation*}
  \begin{array}{c|c|c}
    & \text{half-line} & \text{truncated to } 0<x<15 \\[2pt] \hline
    \rb{-2}{\text{even}} & 
    \rb{-2}{b_k = k\big(k+\frac{1}{2}\big)} &
    \rb{-2}{b_k = 14.0625 + O(1/k^2)} \\[4pt] \hline
    \rb{-2}{\text{full}} &
    \rb{-2}{b_k = \frac{1}{6}(k+1) + \frac{17}{72}(k+1)^{-1} + O(1/k^3)} &
    \rb{-2}{b_k = 14.0625 + O(1/k^2)}
  \end{array}.
\end{equation*}
Since $\|p_j\|^2 = c_j = c_0\prod_{k=1}^j b_j$, we see that the monic
polynomial norms grow super-exponentially on the half-line and
exponentially on the truncated domain.  If $n$ is not too large (say
$n<400$), this alleviates the need to use special floating point
numbers to guard against overflow and underflow.  For larger $n$,
there is little advantage in this respect.  Third, orthogonal
polynomials on the truncated domain are more efficient at representing
functions supported near the origin since their zeros remain confined
to the truncated interval, and therefore can resolve more features of
the solution with fewer basis functions.  On the other
hand, an obvious drawback of truncating the domain is
that the solution and its gradient must remain essentially zero at the
right endpoint to remain a good approximation of the solution on the
half-line.

\begin{figure}
\begin{center}
\includegraphics[width=\linewidth,trim=0 20 0 20]{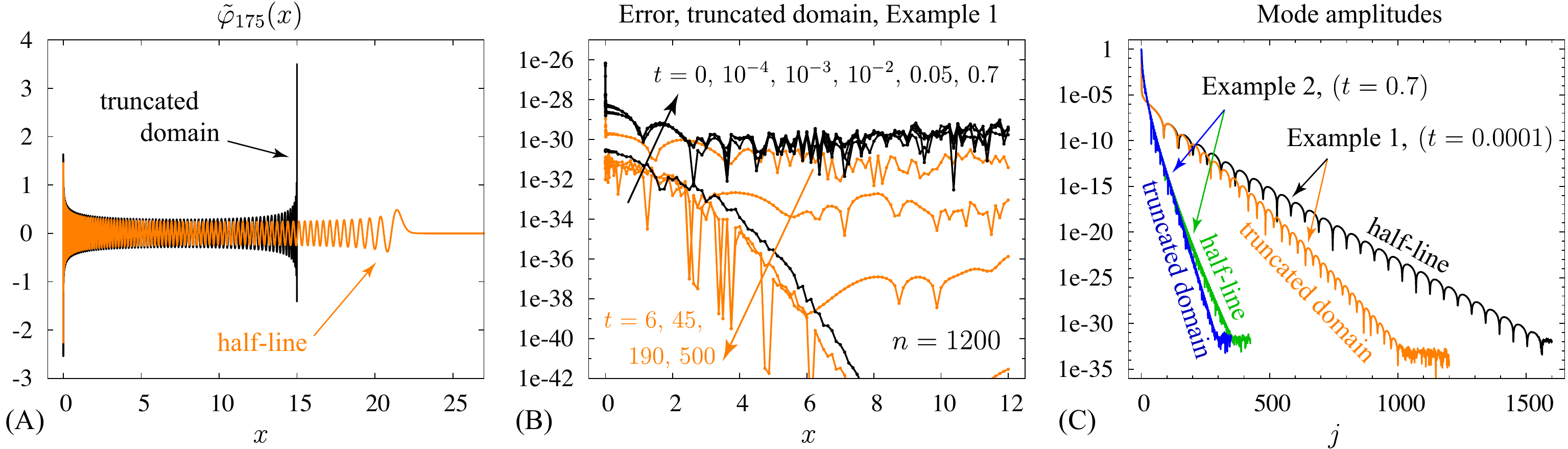}
\end{center}
\caption{\label{fig:trunc1} Plots showing the use of orthogonal
  polynomials on a truncated domain to solve Examples 1 and 2.  The
  error in (B) is nearly identical to the half-line case of
  Figure~\ref{fig:err:compare}(C), but only required 1200 (instead of
  1600) modes. In (C), the modes decay faster on the truncated domain
  because the oscillations in the basis functions are more localized
  near the origin, as shown in (A).  }
\end{figure}

In Figure~\ref{fig:trunc1}, we illustrate these issues in the context
of Examples~1 and~2, where $u(x,0)=x^k$, $k=1$ or $k=2$.  Panel A
compares the 175th normalized basis function, defined in
(\ref{eq:tilde:phi}), on the half-line and truncated domain.  Both
functions have 175 zeros, but they are more spread out in the
half-line case.  Panel B repeats the calculation of
Figure~\ref{fig:err:compare}, showing the difference between the
solution obtained from the projected dynamics (this time on a
truncated domain) to the spectral transform solution (plotted in
Figure~\ref{fig:err:compare}(A,D)).  The errors are essentially the
same as for the half-line (ranging from $10^{-26}$ near $x=0$ to
$10^{-29}$ for larger $x$), but only 1200 modes were needed to reach
roundoff accuracy instead of 1600.  As before, with this choice of
$n$, the projected dynamics solution for Example 1 is correct at $t=0$
and for $t>10^{-4}$, but not at intermediate times $0<t<10^{-4}$.
Similar results, valid for $t\ge0$, were obtained for Example 2.
Panel C shows the mode amplitudes $|\alpha_j(t)|$ at $t=10^{-4}$ in
Example 1 and $t=0.7$ for Example 2.  Fewer modes are needed on the
truncated domain since the zeros remain confined to $0<x<15$ in that
case.

\begin{figure}
\begin{center}
\includegraphics[width=\linewidth,trim=0 20 0 0]{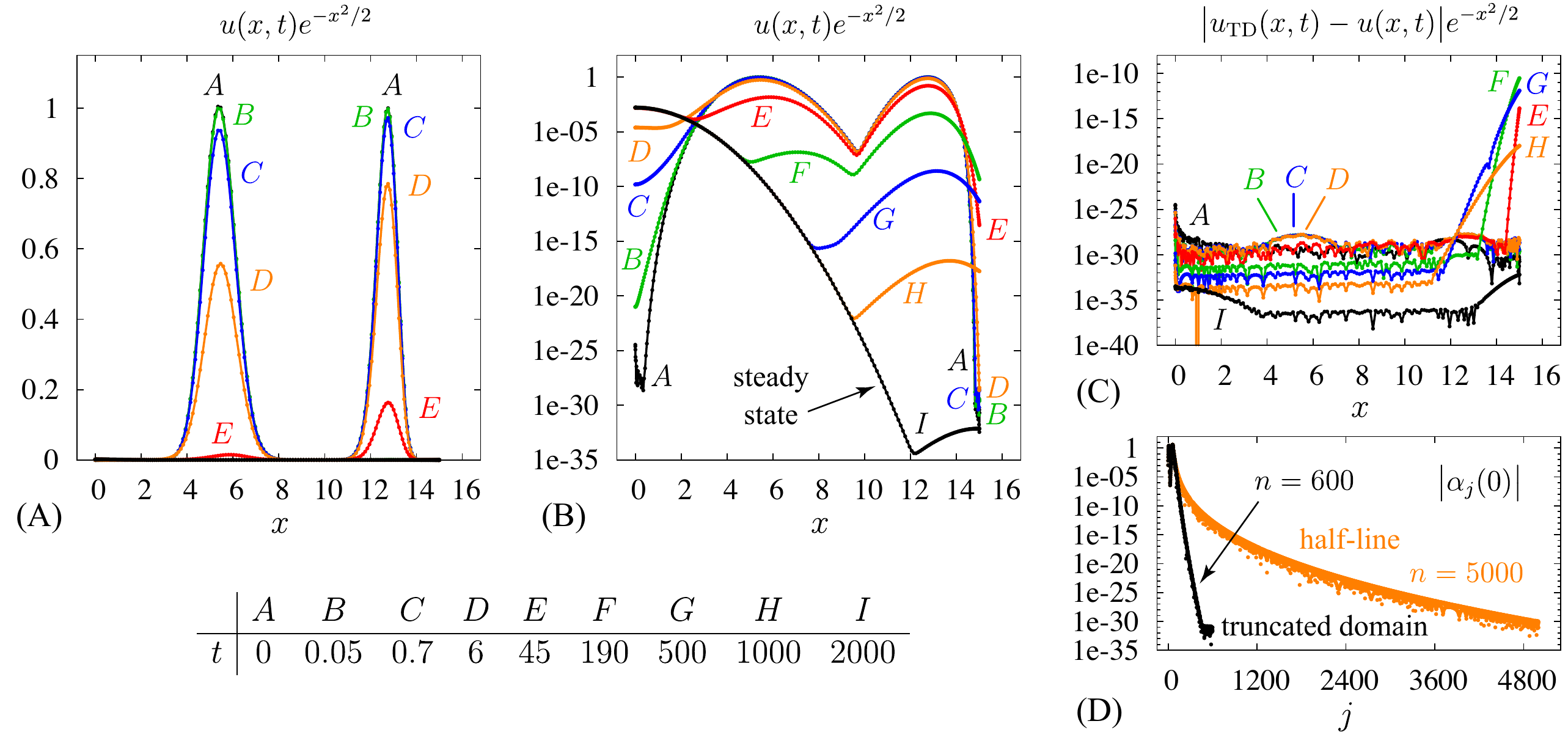}
\end{center}
\caption{\label{fig:trunc2} Evolution of two-hump initial condition
that is initially contained in $0<x<15$ but eventually spreads out
past the right endpoint. In panel (C), $u_\text{TD}(x,t)$ was computed on a
truncated domain while $u(x,t)$ was computed on the half-line. The
modes decay faster in (D) for the truncated domain due to clustering
of the zeros of $p_n(x)$ near $x=15$.}
\end{figure}

Next we look for an initial condition that is initially confined to
$0<x<15$ but spreads out past the right endpoint. We tried a number of
formulas and settled on a two-hump initial distribution of the form
\begin{equation}\label{eq:ex:15}
  u(x,0) = \begin{cases} \left[ \frac{5}{2} \left(\frac{x}{3.25}\right)^{30} +
    \frac{155}{64} \left(\frac{x}{8.25}\right)^{200} \right] \exp\left( -\frac{15}{15-x} \right),
  & 0<x<15, \\
  0, & \text{otherwise}.
  \end{cases}
\end{equation}
The results are shown in Figure~\ref{fig:trunc2}.  Panels A and B show
the solution at the times $t=0,0.05,0.7,6,45,190,500,1000,2000$,
computed on the half-line, on a
linear and log scale, respectively.  Up until $t=6$, the solution
remains confined to $0<x<15$. But then from $6<t<1000$, the solution
is not negligible (in quadruple-precision) at the right endpoint.
Panel C shows the difference between the truncated domain solution and
the half-line solution. As expected, they agree to roundoff error for
$0<t<6$, but then begin to differ near the right end of the domain due
to an incorrect assumption that $u_x=0$ at $x=15$ in the truncated
domain calculation.  Once $t>1000$, the solution has decayed to the
steady-state Maxwellian distribution, and the two methods agree again.
Panel D shows that many fewer modes were needed to resolve the initial
condition on the truncated domain than on the half-line.  This is
because the nodes cluster at $x=0$ and $x=15$ for the truncated domain
calculation, but only at $x=0$ for the half-line calculation, and this
initial condition varies rapidly near $x=15$ due to the factor of
$\exp[15/(15-x)]$ in (\ref{eq:ex:15}).  Even in this example, where
the initial distribution $u(x,t)e^{-x^2/2}$ has a peak of order 1 near
$x=12.75$, the error in truncating the domain to $0<x<15$ was never
larger than $10^{-10}$. Thus, we expect that in practice it is safe to
truncate the domain as long as the initial condition and any sources
are fully supported inside the truncated region.

\begin{figure}
\begin{center}
\includegraphics[width=.95\linewidth,trim=0 20 0 20]{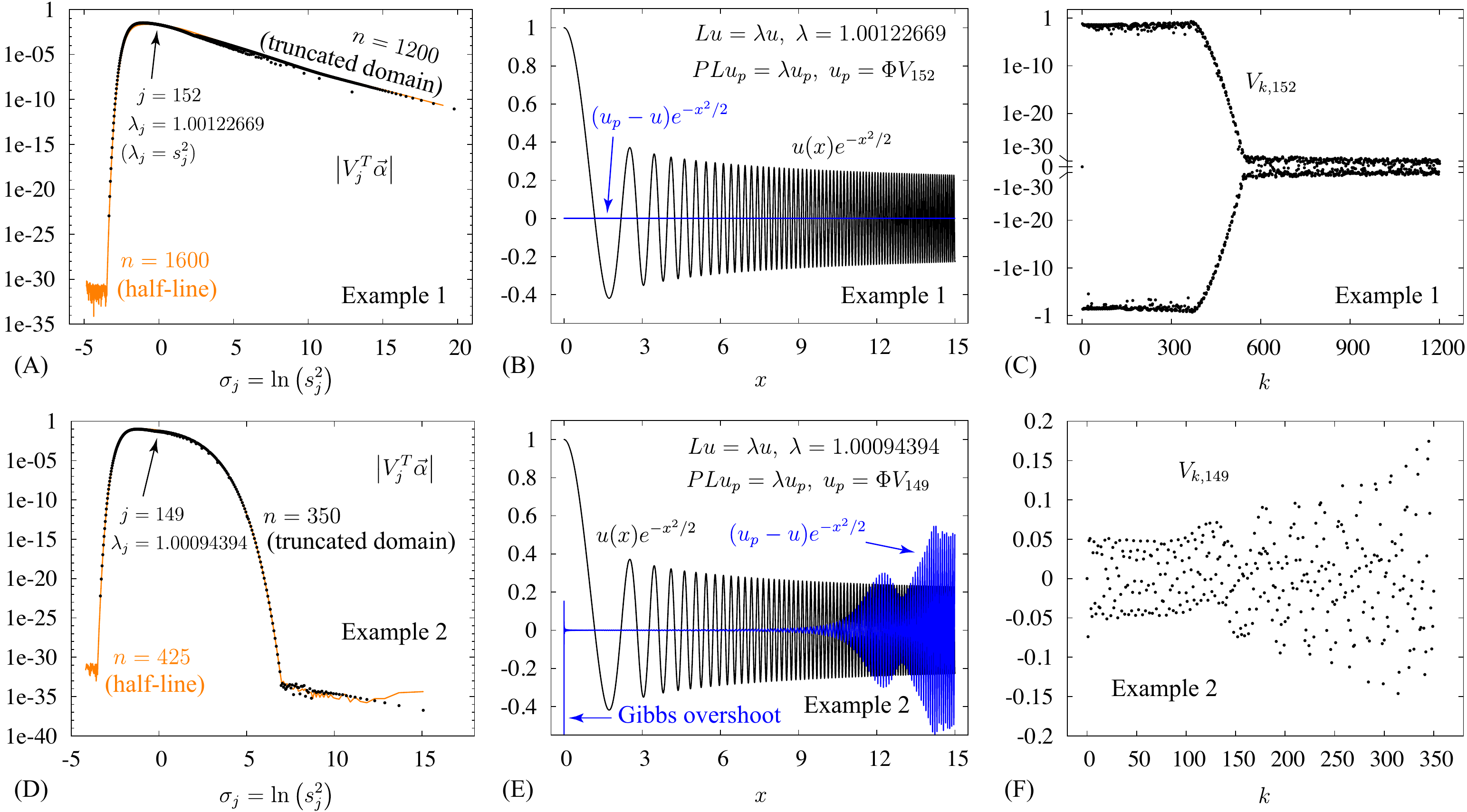}
\end{center}
\caption{\label{fig:trunc3} Comparison of mode amplitudes on a
  half-line and truncated domain (A,D) and plots of resolved (B) and
  unresolved (E) eigenfunctions of $PL$ and $L$ on $0<x<15$. As
  before, $V_j$ is the $j$th column of the matrix $V$ in the singular
  value decomposition $R=USV^T$, and $\Phi V_j$ is the $j$th
  eigenfunction of $PL$, where $\Phi=(\varphi_0,\dots,\varphi_{n-1})$
  are the orthogonal polynomials.}
\end{figure}

As a final remark, we note that the eigenfunctions of $PL$ converge to
eigenfunctions of $L$ when the domain is truncated, unlike the
half-line case.  This is illustrated in Figure~\ref{fig:trunc3}.
Panels A and D show the eigenmode amplitudes for Examples 1 and 2 on a
half-line and truncated domain.  The main difference is that on a
truncated domain, further mesh refinement will not increase the
density of eigenvalues of $PL$ at the left end of the spectrum since
these are converged eigenvalues of the continuous problem. By
contrast, on the half-line, the eigenvalues of $PL$ become more
densely spaced as the mesh is refined to better approximate the
continuous spectrum of $L$.  In Figure~\ref{fig:trunc1}, we saw that
1200 modes was sufficient to resolve the solution in Example 1 and 350
modes was sufficient for Example 2.  In panels (B) and (E) here, we
plot the eigenfunction $u_p$ of $PL$ with eigenvalue $\lambda$
closest to~1, re-scaled to agree as closely as possible with $u$, the
solution of $Lu=\lambda u$ with the same $\lambda$ and
satisfying $u(0)=1$.  Note that 1200 modes yields perfect
agreement between $u$ and $u_p$ (to 26 digits) while 350 modes yields
large discrepancies: a Gibbs overshoot occurs near $x=0$, and a
beat pattern appears for large $x$ where the two functions fall
out of phase.  The minimum amplitude of oscillation of $u_p-u$ in
panel E is roughly $10^{-3}$ when $0.8<x<1.5$, which is 23 orders of
magnitude larger than in panel B.  Panels C and F show why this
occurs: the solution $u$ of $Lu=\lambda u$ with $\lambda\approx 1$
requires about 500 orthogonal polynomial basis functions to be
represented on $0<x<15$.  When 1200 basis functions are used, the
modes decay rapidly between $k=400$ and $k=500$ and the remaining
modes are zero up to roundoff error.  When 350 modes are used, as in
Example 2, the eigenfunction of $PL$ does not agree closely with that
of $L$ due to lack of resolution.  However, just as in the half-line
case, where none of the eigenfunctions of $PL$ agree closely with
solutions of $Lu=\lambda u$, the projected dynamics is accurate to
roundoff error in both examples --- it is not necessary to resolve all
the active eigenfunctions to accurately represent the solution of the
PDE on a truncated domain either.

\section{Conclusion}
\label{sec:conclusion}
Shizgal \cite{shizgal} showed that a new class of non-classical
polynomials could be used very effectively for numerical quadratures
of the collision operator in kinetic simulations. More recently,
Landreman and Ernst \cite{landreman} showed that the same polynomials
could also be much more accurate than other schemes for the
discretization of the speed variable in steady-state kinetic
calculations involving the collision operator. In this article, we
have demonstrated that the new polynomials are also useful for
time-dependent problems by examining the one-dimensional relaxation of a distribution function to a
Maxwell-Boltzmann distribution via energy diffusion.  We found that
the new polynomials are effective at representing the solution of the
partial differential equation for a wide class of initial conditions,
and can be more accurate than generalized Hermite (``even'')
polynomials by many orders of magnitude for the same computational
work. This was seen e.g.~in Fig.~\ref{fig:evol2}(D,E), where 350 modes
are sufficient to reach errors of $10^{-30}$ with the new polynomials
and $10^{-20}$ for the even polynomials, and
Fig.~\ref{fig:err:compare}(C,F), where 1600 modes with the new
polynomials are $10^{20}$ times more accurate than 5000 classical
modes. 

Given the Sturm-Liouville structure that we associated with the
problem, the polynomials defined with integration weight
$\rho(x)=x^{2}e^{-x^2}$ are the most natural to use, and the most
accurate for most parts of the computations. As discussed in \ref{sec:float},
the polynomials defined
with integration weight $\rho(x)=e^{-x^2}$ (on the half-line), chosen for most
computations in \cite{landreman}, also give satisfying results.

It is often the case, given the size of the numerical simulations,
that one can only afford very coarse grids for a given variable.  Our
analysis at low resolution, in Figures \ref{fig:smallGrid} and
\ref{fig:smallGridEntropy}, shows that the full polynomials are more
accurate than even polynomials as well as the discretization scheme
used in popular plasma kinetic codes, often by several orders of
magnitude. These results, together with exact mass
conservation at all times suggests that the new polynomials could be
an attractive alternative to the finite difference schemes currently
in use in state-of-the-art plasma microturbulence codes
\cite{candy1,barnes1}. However, for the new polynomials to make a
truly compelling case, at least two questions must be answered.
First, the stiffness matrix $K$ is dense, whereas finite difference
matrices are sparse. Does the fact that the polynomials yield accurate
results on very coarse grids compensate this disadvantage?  Are there
formulations based on these polynomials that can avoid operations on
dense matrices?  Curiously, when $\nu=0$, we find that $M$ is
pentadiagonal and the entries $K_{ij}$ decay exponentially as $|i-j|$
increases, but this still leads to a much wider band of nonzero
entries centered about the diagonal than in a finite-difference
approach.  Second, what is the best way to incorporate these
polynomials in time-dependent simulations with more complete collision
operators?  Exponential time-differencing schemes \cite{kassam} and
implicit-explicit Runge-Kutta methods \cite{carpenter} appear
promising, but this is the subject of ongoing research, with results
to be reported at a later date.

\section{Acknowledgments}
J.W.~was supported in part by the U.S. Department of Energy, Office of
Science, Office of Advanced Scientific Computing Research, Applied
Mathematics program under contract number~DE-AC02-05CH11231, and by
the National Science Foundation under Grant No.~DMS-0955078.
A.J.C.~was supported by the U.S. Department of Energy, Office of
Science, Fusion Energy Sciences under Award No.~DE-FG02-86ER53223 and
DE-SC0012398. M.L.~was supported in part by the U.S. Department of
Energy, Office of Science, Office of Fusion Energy Science, under
award numbers DE-FG02-93ER54197 and DE-FC02-08ER54964.

\appendix

\section{Effect of the choice of $\nu$ in floating-point arithmetic}
\label{sec:float}

In exact arithmetic, the projected dynamics onto
$\V=\opn{span}\{\varphi_0,\dots,\varphi_{n-1}\}$ is identical for any
choice of $\nu$ in (\ref{eq:rho:plus}) since the subspace $\V$
consists of polynomials of degree less than or equal to $n$,
regardless of the weight function.  However, changing $\nu$ can have a
large effect in the presence of roundoff error.  To understand this,
first observe that using a monomial basis for $\V$ (instead of
orthogonal polynomials) would lead to numerical difficulties as the
basis functions become nearly linearly dependent in
$\H=L^2(\mathbb{R}_+;w\,dx)$.  Using orthogonal polynomials with
weight exponent $\nu\ne2$ has the potential to cause similar
difficulties.

In exact arithmetic, the mass matrix $M=R_1^TR_1$ in (\ref{eq:MK}) is
the identity when $\nu=2$ since the $\varphi_j$ in (\ref{eq:varphi})
are orthonormal in $\H$.  For any other choice of $\nu$, the Cholesky
factor $R_1$ gives the change of basis to the $\nu=2$ case:
\begin{equation*}
  \Phi^{(\nu)} =
  \big(\varphi_0^{(\nu)},\dots,\varphi_{n-1}^{(\nu)}\big)
  = \Phi^{(2)}R_1^{(\nu)}, \qquad
  M^{(\nu)} = \int \Phi^{(\nu)\,T}\Phi^{(\nu)}w\,dx
  = R_1^{(\nu)\,T}R_1^{(\nu)}.
\end{equation*}
Similarly, $K^{(\nu)} = R_1^{(\nu)\,T}K^{(2)}R_1^{(\nu)}$ so that $R =
R_2^{(\nu)}(R_1^{(\nu)})^{-1} = R_2^{(2)}$ in (\ref{eq:R12}) is
independent of $\nu$ (up to roundoff errors), and is largely
unaffected by orthogonality drift when $\nu=2$ as long as
$R_1^{(\nu)}$ is actually computed rather than assumed to be the
identity matrix.  From (\ref{eq:exp:MK}), the solution of the
projected dynamics is
\begin{equation}\label{eq:proj:soln2}
  u_p(x,t) = \big[ \underbrace{\Phi^{(\nu)}(x)(R_1^{(\nu)})^{-1}}_{A}\big]
  \big[ V e^{-{S^2}t} V^T\big] \big[
  \underbrace{R_1^{(\nu)}\vec\alpha^{(\nu)}}_{B}\big]
\end{equation}
where
$A=\Phi^{(2)}(x)$, $B = \vec\alpha^{(2)} = (R_1^{(\nu)})^{-T}\vec\beta^{(\nu)}$,
$\vec\beta^{(\nu)}=M^{(\nu)}\vec\alpha^{(\nu)}$, and
\begin{equation}\label{eq:alpha:beta}
\alpha^{(\nu)}_i = \int u_p(x,0)\,
\overline{\varphi^{(\nu)}_i(x)}\,x^\nu e^{-x^2}dx, \qquad
\beta^{(\nu)}_i = \int u_p(x,0)\,
\overline{\varphi^{(\nu)}_i(x)}\,x^2 e^{-x^2}dx.
\end{equation}
Thus, if $\nu\ne2$, the numerical algorithm performs a change of basis
to the $\nu=2$ case as an intermediate step.

\begin{figure}[t]
\begin{center}
\includegraphics[width=\linewidth,trim=0 20 0 20]{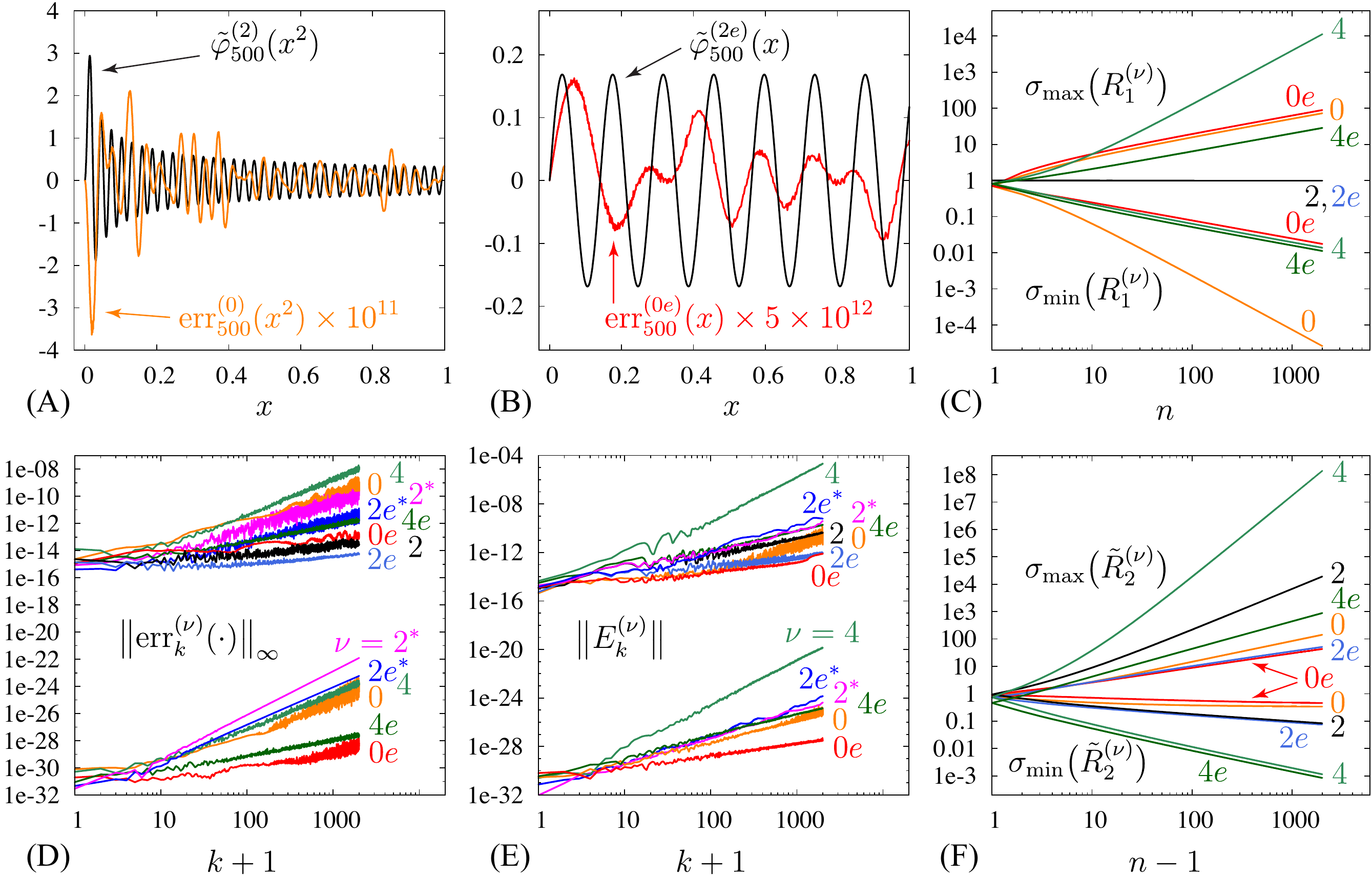}
\end{center}
\caption{\label{fig:phi2err} Comparison of roundoff errors for
  $\nu=0,0e,2,2^*,2e,2e^*,4,4e$. Here $e$ stands for ``even,'' $*$
  stands for ``assuming $R_1=I$,'' and errors are measured against the
  $\nu=2$ or $\nu=2e$ quadruple-precision results.  (A) and (B) show
  the error in using $\tilde\Phi^{(\nu)}(x)(R_1^{(\nu)})^{-1}$ with
  $\nu=0$ to compute $\tilde\varphi_{500}^{(2)}(x)$.  (D) gives the
  sup-norm of the error in computing $\tilde\varphi_{k}^{(2)}(x)$ in
  this way.  (E) shows the error in the $k$th column of
  $\opn{pinv}(R)$, as defined in (\ref{eq:E:nu:def}). (C) and (F) give
  the largest and smallest singular values of $R_1$ and $\tilde R_2$,
  which help explain the results in (A), (B), (D), (E).}
\end{figure}

Figure~\ref{fig:phi2err} shows the sources and effect of roundoff
error in evaluating the solution $u_p(x,t)$ via
(\ref{eq:proj:soln2}).  Panels (A) and (B) show the normalized basis
functions $\tilde\varphi_{500}^{(2)}(x)$ and
$\tilde\varphi_{500}^{(2e)}(x)$ in (\ref{eq:tilde:phi}) and the error
in computing them in double-precision with $\nu=0$ or $\nu=0e$, where
``$e$'' stands for even:
\begin{equation}\label{eq:err:def}
  \text{err}_{k}^{(\nu)}(x) = \opn{col}_k\big[
  \tilde\Phi^{(\nu)}(x)\big(R_1^{(\nu)}\big)^{-1}\big] -
  \tilde\varphi_k^{(2\text{ or }2e)}(x).
\end{equation}
Here $\opn{col}_k X$ is the $k$th column of $X$.
The second term on the right is the ``exact'' solution, which was
computed with $\nu=2$ in quadruple-precision. The relative error for
$k=500$ is
about 60 times smaller in the even case (B) than the full case (A).
This is because the orthogonal polynomials are more oscillatory near
$x=0$ in the full case, where $x=0$ is a true integration boundary.  A
similar thing occurs in Chebyshev and Legendre polynomials, which are
more oscillatory near $x=\pm1$ than $x=0$. We plotted
$\tilde\varphi_{500}^{(2)}(x^2)$ in (A) to obtain more uniform
oscillations for visualization.

In both (A) and (B), we note that the error has more structure than
might be expected from roundoff errors alone. This is because
$(R_1^{(\nu)})^{-1}$ is present in (\ref{eq:err:def}), which amplifies
errors along singular vectors corresponding to the smallest singular
values more than in other directions. (C) shows the largest and
smallest singular values of $R_1^{(\nu)}$ as a function of $n=\dim
\V$.  Since the matrices $R_1^{(\nu)}$ are nested as $n$ increases,
$\sigma_\text{max}$ is an increasing function of $n$ while
$\sigma_\text{min}$ is decreasing.  The condition number of
$R_1^{(\nu)}$ is the ratio $\sigma_\text{max}/\sigma_\text{min}$. It
is 1 for $\nu=2$, grows slowly for $\nu=0e$ and $\nu=4e$, and grows
faster for $\nu=0$ and $\nu=4$.  The condition number of $M$ is the
square of the condition number of $R_1$. Thus, roundoff errors are
reduced by computing $R_1$ directly rather than from $M$.

Panel (D) shows the max-norm of $\text{err}^{(\nu)}_k(x)$ as a
function of $k$ for $\nu=0,2,4$ in the even and full cases. For
example, the amplitudes of the largest peaks in the error curves in
(A) and (B) are plotted at $k=500$ in (D) on the double-precision
curves labeled $0$ and $0e$, respectively. For small values of $k$,
all the values of $\nu$ that we tested yield reasonably accurate
results. However, as $k$ increases, the growth in the condition number
of $R_1^{(\nu)}$ in (C) leads to significant loss of accuracy in
$\tilde\varphi_k^{(2)}$ when computed via
$\tilde\Phi^{(\nu)}(R_1^{(\nu)})^{-1}$ for $\nu\ne2$.  In particular,
roundoff error is amplified by more than 5 orders of magnitude with
$\nu=4$ and $k=2000$.  The $\nu=2$ curves (full and even) are missing
from the quadruple-precision results as they are treated as ``exact''
solutions in (\ref{eq:err:def}).  Note that even for $\nu=2$ we
compute $R_1$ and apply its inverse to $\Phi^{(2)}$ to correct for the
slight loss of orthogonality that occurs when computing the basis
functions by the three-term recurrence (\ref{eq:3:term:brief}) or
(\ref{eq:even:recur}).  In other words, the second term in
(\ref{eq:err:def}) is computed as
\begin{equation}
  \tilde{\varphi}_k^{(2\text{ or }2e)}(x) := \opn{col}_k\big[
  \tilde\Phi^{(2\text{ or }2e)}(x)\big(R_1^{(2\text{ or }2e)}
  \big)^{-1}\big].
\end{equation}
The curves labeled $2^*$ and $2e^*$ give uncorrected results in which
$R_1$ is assumed equal to the identity and dropped from the first term
in (\ref{eq:err:def}).  These errors are comparable to using $\nu=0$
or $\nu=4$, where $R_1$ is required.

Panel (E) shows the error in computing $\opn{pinv}(R)$ in (\ref{eq:pinv}).
More precisely, we plot
\begin{equation}\label{eq:E:nu:def}
  E^{(\nu)}_k = \opn{col}_k\left[ \tilde R^{(\nu)}_1(\tilde R^{(\nu)}_2)^{-1} -
    \tilde R^{(2 \text{ or } 2e)}_1(\tilde R^{(2 \text{ or } 2e)}_2)^{-1} \right],
\end{equation}
where $\tilde R_j$ is obtained from $R_j$ by deleting the zeroth row
and column.  The second term on the right is computed in quadruple
precision, and treated as ``exact.''  Since $R_1$ and $R_2$ are
upper-triangular, the $k$th column of $E^{(\nu)}$ is independent of
$n=\dim \V$ for $n>k$.  Note that the double-precision curves
labeled $\nu=0$ and $\nu=0e$ are more accurate than the $\nu=2$ and
$\nu=2e$ curves, respectively. This is because we compute the SVD of
$\opn{pinv}(R)$, which involves $R_2^{-1}$ rather than $R_1^{-1}$, and
$R_2^{(0)}$ is better conditioned than $R_2^{(2)}$, as shown in panel
(F).  The reason is that the orthogonal polynomials for $\nu=0$ are
less oscillatory near $x=0$ than for $\nu=2$ (since $\rho(x)$ vanishes
at $x=0$ in the latter case), and $R_2$ involves derivatives of these
orthogonal polynomials. Thus, a better ``exact'' solution would be
$\nu=0$ instead of $\nu=2$ in (\ref{eq:E:nu:def}). When this is done,
the only visible effect on the plot in panel (E) is that the two most
accurate quadruple-precision curves should be re-labeled (0 to 2, $0e$
to $2e$) to better account for the primary source of error in
(\ref{eq:E:nu:def}).

As mentioned previously, we compute the SVD of $\opn{pinv}(R)$ rather
than of $R$ because $\|R\|>\|\opn{pinv}(R)\|$.  This may be
seen in panel (F), where $\sigma_\text{max}(\tilde R_2^{(2)})$ grows
faster than $\sigma_\text{min}(\tilde R_2^{(2)})$ decays. (Recall that
$R=R_2^{(2)}$ since $R_1^{(2)}=I$.)  Also, the largest singular values
of a matrix are computed with the most relative accuracy, and the
largest singular values of $\opn{pinv}(R)$ are the ones that matter
most in accurately representing solutions of the PDE (\ref{eq:PDE}).

In summary, the most accurate computation of $u_p(x,t)$ in
(\ref{eq:proj:soln2}) in floating point arithmetic would involve
computing $A=\Phi^{(\nu)}(x)(R_1^{(\nu)})^{-1}$ with $\nu=2$ (without
assuming $R_1^{(2)}=I$), and computing $V$ and $S$ from $\tilde
R_1^{(\nu)}(\tilde R_2^{(\nu)})^{-1}$ with $\nu=0$. In practice we set
$\nu=2$ when computing $V$ and $S$ as well since the improvement in
switching to $\nu=0$ is small. The last term in (\ref{eq:proj:soln2}),
$B=\vec\alpha^{(2)}$, can often be computed analytically.  If not, then it is
most accurately computed as $B=(R_1^{(\nu)})^{-T}\vec\beta^{(\nu)}$ with
$\nu=2$, including $R_1$ as before to account for the slight loss of
orthogonality in the basis functions.

\vspace*{10pt}
\bibliographystyle{elsarticle-num}


\end{document}